\documentclass{gtart_h}  

\def\ifplaintex{\expandafter\ifx\csname documentclass\endcsname\relax}

\def\ifplaintex{\expandafter\ifx\csname documentclass\endcsname\relax}

\ifplaintex 
\else
\fi

\expandafter\ifx\csname epsfbox\endcsname\relax\input epsf\fi

\def\gt{{\mathsurround=0pt\it $\cal G\mskip-2mu$eometry \&\ 
$\cal T\!\!$opology}}        

\def\gtp{{\mathsurround=0pt\it $\cal G\mskip-2mu$eometry \&\ 
$\cal T\!\!$opology $\cal P\!$ublications}}  

\def\lognumber#1{\def\thelognumber{#1}}
\def\volumenumber#1{\def\thevolumenumber{#1}}
\def\papernumber#1{\def\thepapernumber{#1}}
\def\volumeyear#1{\def\thevolumeyear{#1}}

\def\pagenumbers#1#2{\def\startpage{#1}\def\finishpage{#2}}
\def\published#1{\def\publishdate{#1}}
\def\proposed#1{\def\theproposer{#1}}
\def\seconded#1{\def\theseconders{#1}}
\def\received#1{\def\receiveddate{#1}}
\def\revised#1{\def\reviseddate{#1}}
\def\accepted#1{\def\accepteddate{#1}}
\def\asciititle#1{\def\theasciititle{#1}}

\def\asciiemail#1{\def\theasciiemail{#1}}

\long\def\asciiabstract#1{\long\def\theasciiabstract{#1}}
\def\asciikeywords#1{\def\theasciikeywords{#1}}

\let\\\par\let\thelognumber\relax
\let\thevolumenumber\relax\let\thepapernumber\relax
\let\thevolumeyear\relax\let\thesamplenumber\relax\let\startpage\relax
\let\finishpage\relax\let\publishdate\relax\let\receiveddate\relax
\let\reviseddate\relax\let\accepteddate\relax\let\theasciititle\relax
\let\theasciiauthors\relax
\let\theasciiabstract\relax\let\theasciikeywords\relax
\let\theasciiemail\relax\let\theshortauthors\relax\let\theshorttitle\relax

\long\def\maketitlep{   

\count0=\startpage

\gt\hfill      
\hbox to 77pt{\vbox to 0pt{\vglue -15pt\epsfbox{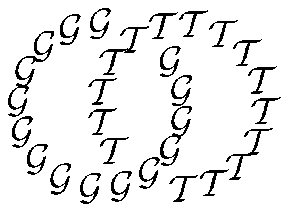}\vss}\hss}
\break
{\small\ifx\thesamplenumber\relax 
Volume \else Sample
\fi\thevolumenumber\ (\thevolumeyear)
\startpage--\finishpage\nl
Published: \publishdate}
\vglue 0.5truein plus 0.4fil minus 0.1truein

{\parskip=0pt\leftskip 0pt plus 1fil\def\\{\par\smallskip}{\ifplaintex\large
\else\Large\fi\bf\thetitle}\par\medskip}   

\vglue 0pt plus 0.1fil 

{\parskip=0pt\leftskip 0pt plus 1fil\def\\{\par}{\sc\theauthors}
\par\medskip}

\vglue 0pt plus 0.1fil 

{\small\parskip=0pt\let\newline\\
{\leftskip 0pt plus 1fil\def\\{\par}{\sl\theaddress}\par}
\expandafter\ifx\theemail\relax    
\relax\else\vglue 5pt plus 0.02fil minus 2pt\def\\{\stdspace{\rm 
and}\stdspace} 
\cl{Email:\stdspace\tt\theemail}\fi
\ifx\theurl\relax                  
\relax\else\vglue 5pt plus 0.02fil minus 2pt\def\\{\stdspace{\rm 
and}\stdspace}
\cl{URL:\stdspace\tt\theurl}\fi\par}

\vglue 7pt plus 0.3fil minus 3pt

{\bf Abstract}
\vglue 5pt plus 0.1fil minus 2pt

\theabstract

\vglue 7pt plus 0.3fil minus 3pt

{\bf AMS Classification numbers}\quad Primary:\quad \theprimaryclass

Secondary:\quad \thesecondaryclass

\vglue 5pt plus 0.3fil minus 2pt

{\bf Keywords:}\quad \thekeywords

\vglue 10pt plus 0.5fil minus 5pt

{\small  Proposed: \theproposer\hfill Received: \receiveddate\nl
Seconded: \theseconders\hfill 
\ifx\reviseddate\relax                         
Accepted: \accepteddate                        
\else
Revised: \reviseddate                          
\fi}
\eject
}       

\font\phead=cmsl9 scaled 950
\font=cmsl9 scaled 1050
\font\pnum=cmbx10 scaled 913
\font\lnum=cmbx10 
\font\pfoot=cmsl9 scaled 950
\font=cmsl9 scaled 1050
\ifplaintex
\headline{\vbox to 0pt{\vskip -4.5mm\line{\small\phead\ifnum
\count0=\startpage ISSN 1364-0380 (on line)
1465-3060 (printed) \hfill {\pnum\folio}\else\ifodd\count0\def\\{ }%
\ifx\theshorttitle\relax\thetitle\else\theshorttitle\fi\hfill{\pnum\folio}
\else\def\\{ and }{\pnum\folio}\hfill\ifx\theshortauthors\relax\theauthors
\else\theshortauthors\fi\fi\fi}\vss}}
\footline{\vbox to 0pt{\vglue 0mm\line{\small\pfoot\ifnum\count0=\startpage
\copyright\ \gtp\hfill\else
\gt, Volume \thevolumenumber\ (\thevolumeyear)\hfill\fi}\vss
}}
\else
\makeatletter
\def\@oddhead{{\small\lhead\ifnum\count0=\startpage ISSN 1364-0380 (on line)
1465-3060 (printed) \hfill {\lnum\number\count0}\else\ifodd\count0
\def\\{ }\ifx\theshorttitle\relax \thetitle \else\theshorttitle\fi\hfill
{\lnum\number\count0}\else\def\\{ and }{\lnum\number\count0}
\hfill\ifx\theshortauthors\relax 
\theauthors\else\theshortauthors\fi\fi\fi}}\def\@evenhead{\@oddhead}
\def\@oddfoot{\small\lfoot\ifnum\count0=\startpage\copyright\ \gtp\hfill\else
\gt, Volume \thevolumenumber\ (\thevolumeyear)\hfill\fi}
\def\@evenfoot{\@oddfoot}
\makeatother
\fi

\newwrite\gtoutfile
\long\gdef\makeheadfile{  
{\def\\{, }\def\s{ }
\immediate\openout\gtoutfile head.xxx
\immediate\write\gtoutfile{Proxy-for: \ifx\theasciiauthors\relax
\theauthors\else\theasciiauthors\fi\s<\ifx\theasciiemail\relax\theemail\else\theasciiemail\fi>}
\immediate\write\gtoutfile{\noexpand\\}
\immediate\write\gtoutfile{Authors: \ifx\theasciiauthors\relax
\theauthors\else\theasciiauthors\fi}
{\def\\{ }\immediate\write\gtoutfile{Title: \ifx\theasciititle\relax
\thetitle\else\theasciititle\fi}}
\immediate\write\gtoutfile{Subj-class: GT or SG or MG etc}
\immediate\write\gtoutfile{MSC-class: \theprimaryclass\ifx\thesecondaryclass\relax\else, \thesecondaryclass\fi}
\immediate\write\gtoutfile{Journal-ref: Geom. Topol. \thevolumenumber
(\thevolumeyear) \startpage-\finishpage}
\immediate\write\gtoutfile{Comments: Published by Geometry and Topology at}
\immediate\write\gtoutfile{\s\s http://www.maths.warwick.ac.uk/gt/GTVol\thevolumenumber/paper\thepapernumber.abs.html}
\immediate\write\gtoutfile{\noexpand\\}
\immediate\write\gtoutfile{}
\ifx\theasciiabstract\relax
\immediate\write\gtoutfile{\theabstract}\else
\immediate\write\gtoutfile{\theasciiabstract}\fi
\immediate\write\gtoutfile{}
\immediate\write\gtoutfile{\noexpand\\}
\immediate\write\gtoutfile{}
\immediate\closeout\gtoutfile}}  

\def\maketitlepage{\maketitlep\makeheadfile}
\let\maketitle\maketitlepage

\lognumber{347}
\volumenumber{8}\papernumber{11}\volumeyear{2004}
\pagenumbers{475}{509}
\received{29 July 2003}
\revised{17 January 2004}
\published{16 February 2004}
\accepted{9 February 2004}
\proposed{Ronald Fintushel}
\seconded{Ronald Stern, Robion Kirby}

\usepackage{amsmath,amssymb}
\usepackage[all]{xy}

\usepackage[dvips,bookmarks,bookmarksnumbered,pdfstartview=FitBH]{hyperref}

\def\S{Section }
\newtheorem{dummy}{anything}[section] 
\newtheorem{theorem}[dummy]{Theorem}
\newtheorem*{thma}{Theorem A}
\newtheorem*{thmb}{Theorem B}
\newtheorem*{thmc}{Theorem C}
\newtheorem{lemma}[dummy]{Lemma} 
\newtheorem{proposition}[dummy]{Proposition} 
\newtheorem{corollary}[dummy]{Corollary}
 
\theoremstyle{definition}
\newtheorem{definition}[dummy]{Definition}
\newtheorem{example}[dummy]{Example}
\newtheorem{remark}[dummy]{Remark}
\newtheorem*{rem}{Remark}
\newtheorem*{ack}{Acknowledgement}

\newbox\nbox
\newcommand{\bbar}[1]
{\setbox\nbox = \hbox{$#1$}
\kern .2\wd\nbox{\overline
{\hbox to .7\wd\nbox{\vphantom {$#1$}\hss}}}
\kern .1\wd\nbox\kern -\wd\nbox\box\nbox}
\newcommand{\wbar}[1]{\overset{\, \hrulefill\, }{#1}}

\newcommand{\vv}{\: | \:}
\DeclareMathOperator{\Ker}{Ker} 
\DeclareMathOperator{\Coker}{Coker}
\DeclareMathOperator{\Fix}{Fix}
\DeclareMathOperator{\SH}{\#}
\DeclareMathOperator{\ind}{Ind}
\DeclareMathOperator{\ad}{ad}
\newcommand{\OmegaP}
{\Omega^*_{-}(\ad P)}
\DeclareMathOperator{\Sing}{Sing}

\newcommand{\G}{\pi}
\renewcommand{\H}{\pi}

\newcommand{\bd}{\partial\,}
\newcommand
{\mmatrix}[4]{\left (\vcenter
{\xymatrix@C-2pc@R-2pc
{#1&#2\\#3&#4}}\right )}
\newcommand{\hp}{{\hphantom{-}}}
\newcommand{\bT}{\mathbf T}

\newcommand{\bbR}{\mathbb R}
\newcommand{\bbT}{\mathbb T}
\newcommand{\bC}{\mathbb C}

\newcommand{\bR}{\mathbb R}
\newcommand{\bZ}{\mathbf Z}

\newcommand{\shp}{\,\#\,}
\newcommand{\cp}{\bC P^2}
\newcommand{\nCP}{\SH_1^n \bC P^2}
\newcommand{\cA}{\mathcal A}
\newcommand{\cB}{\mathcal B}
\newcommand{\cC}{\mathcal C}
\newcommand{\cD}{\mathcal D}

\newcommand{\cF}{\mathcal F}
\newcommand{\cG}{\mathcal G}

\newcommand{\cM}{\mathcal M}
\newcommand{\Me}{\wbar{\cM^*_{(e)}}}
\newcommand{\cN}{\mathcal N}
\newcommand{\cR}{\mathcal R}
\newcommand{\cS}{\mathcal S}
\newcommand{\cU}{\mathcal U}

\DeclareMathOperator{\lk}{\ell}

\begin{document}
\title[Permutations, isotropy and smooth actions on $4$--manifolds]{Permutations,  isotropy and smooth cyclic\\group actions on definite $4$--manifolds}
\asciititle{Permutations, isotropy and smooth cyclic group actions on definite 4-manifolds}

\authors{Ian Hambleton\\Mihail Tanase}
\address{Department of Mathematics and Statistics\\McMaster 
University, Hamilton, ON L8S 4K1, Canada}
\gtemail{\mailto{ian@math.mcmaster.ca}\qua{\rm and}\qua 
\mailto{tanasem@math.mcmaster.ca}}
\asciiemail{ian@math.mcmaster.ca, tanasem@math.mcmaster.ca}

\begin{abstract}
We use the equivariant Yang--Mills moduli space to investigate the relation between the singular set, isotropy representations at fixed points, and permutation modules realized by the induced action on homology for smooth group actions on certain $4$--manifolds.
\end{abstract}

\asciiabstract{We use the equivariant Yang-Mills moduli space to
investigate the relation between the singular set, isotropy
representations at fixed points, and permutation modules realized by
the induced action on homology for smooth group actions on certain
4-manifolds.}

\asciikeywords{Gauge theory, 4-manifolds, group actions, Yang-Mills,
moduli space}
\maketitlepage

\section*{Introduction}
We are interested in a kind of ``rigidity"  for  finite group actions on certain smooth $4$--manifolds, namely those constructed by connected
sums of  geometric pieces such as algebraic surfaces.  We can ask
how closely a smooth, orientation-preserving,  finite group action on such a connected sum resembles a equivariant connected sum of algebraic actions on
the individual factors.

Following  \cite{HL1},  we consider this question for the simplest case, where
$$X = \nCP = \cp \shp \dots \shp \cp$$ is the connected sum of $n$
copies of the complex projective plane.
 In that paper we restricted ourselves to actions
which induced the \emph{identity} on the homology of $X$,  but
here we remove this assumption. Since the intersection form
of $X$ is the standard definite form $Q_X = \langle 1\rangle\perp
\dots \perp\langle 1\rangle$, we note that its automorphism
group is given by an extension
$$1 \to \{\pm 1\}^n \to \text{Aut} (H_2(X;\bZ), Q_X) \to \Sigma_n\to 1$$
where $\Sigma_n$ denotes the group of permutations of $n$
elements. Therefore, if a finite group $\G$ of odd order acts smoothly on $X$,  then $H_2(X;\bZ)$ is the direct sum of permutation modules
of the form $\bZ[\G/\H_\alpha]$, for various stabilizer subgroups
$\H_\alpha \subseteq \G$.
\begin{rem} Our results (like those of \cite{HL1}) actually hold for any smooth, closed, 
simply-connected $4$--manifold with positive definite intersection form.
This is basically due to the result of S Donaldson \cite{Doa}, 
that
any such $4$--manifold $X$ has the standard positive definite intersection form,
and hence is homotopy equivalent to $\nCP$ for some $n\geq 0$.
With this information, it is easy to see that the arguments go through
without essential changes. 
\end{rem}

We would like to understand how the following three invariants
of such an action $(X, \G)$ are related:

\begin{enumerate}
\renewcommand{\labelenumi}{(\Alph{enumi})}
\item The permutation representation of $\G$ on $H_2(X;\bZ)$.
\item The singular set of the action, meaning the collection
of isotropy subgroups and fixed sets $\Fix(X,\H')$ for $\H' \subset \G$.
\item The tangential isotropy representations $(T_x X, \G_x)$
at all singular points $x \in X$.
\end{enumerate}

In the rest of the paper, we assume that $\G = C_m$ is a finite
cyclic group of \emph{odd} order $m$, acting smoothly on
$X = \nCP$.  
We  obtain many examples of such smooth actions
by starting with linear actions of $\G$ on $\cp$.  These are just the cyclic subgroups
of the algebraic automorphism group $PGL_3(\bC)$, given by sending
a generator 
$$ t\ \mapsto\ \left ( \vcenter{\xymatrix@R-28pt@C-24pt{
1&&\cr&\zeta^a&\cr&&\zeta^b}}\right )$$
where $\zeta = e^{\frac{2\pi i}{m}}$ is a primitive $m^{th}$ root
of unity, and $a$, $b$ are integers such that the greatest common divisor
$(a,b,m) =1$. In this case, $\G$ acts by the identity on homology,
and the singular set always contains the three fixed points $[1,0,0]$,
$[0,1,0]$, and $[0,0,1]$. In addition, there can be up to three invariant
$2$--spheres with various isotropy subgroups, depending on the values of $a$ and $b$. For example, if $(a,b) = (10,3)$ and $m = 105$,
then the action has 5 orbit types (the maximal number for $\cp$).
The tangential representations at the three $\G$--fixed points
are given by the \emph{rotation numbers} $(a,b)$, $(-a, b-a)$
and $(a-b, -b)$ standing for the decomposition of $T_x \cp =
\bR^4 = \bC^2$ into eigenspaces under the action of $t$.
These rotation numbers are well-defined modulo $m$ up to
identifying $(a,b) \equiv (b,a) \equiv (-a,-b)$.

Examples of smooth $\G$--actions on a connected sum $X =\nCP$
with more isotropy groups and various permutation actions on homology are constructed
by a \emph{tree} of equivariant connected sums, where we connect
either at fixed points of two linear actions
or along an orbit of singular points. In order to preserve orientation,
the rotation numbers at the attaching points must be of the
form $(a,b)$ and $(a,-b)$. By this means, we can obtain a large supply
of model actions, for which one should be able to work out the
relation between (A), (B) and (C). 

What can one say about a general smooth action $(X,\G)$ ? It is
not hard to verify that the singular set consists of a configuration
of isolated points and $2$--spheres, as in the linear models (see
\cite{Ed}).
 The main result of  \cite{t}
generalizes \cite[Theorem C]{HL1}:
\begin{thma}
 Let $X$ be a smooth, closed, simply-connected
$4$--manifold with a positive definite intersection form.
 Let $\G$ be a  cyclic group of odd order, acting smoothly on 
on $X$. Then there exists an
equivariant connected sum of linear actions on $\cp$ with the same
isotropy structure, singular set and rotation numbers, and the same
permutation action on $H_2(X;\bZ)$ as for the given action $(X, \G)$.
\end{thma}
This result is proved by using the symmetries of the equivariant Yang--Mills moduli space \cite{HL2} to produce a stratified equivariant cobordism 
between $(X,\G)$ and a connected sum of linear actions, relating
the invariants (A), (B) and (C). The
following result of \cite{t}  holds for smooth, but not for topological actions, by an example
of A Edmonds \cite{Ed1}.
\begin{thmb}
 Let $(X,\G)$ be a smooth action as above, with
discrete singular set. Then the action $(X,\G)$ is semi-free.
\end{thmb}

In principle, it should now be possible to say exactly which permutation
modules can be realized  by smooth actions on $\nCP$ just
by studying the equivariant connected sums. We will show
that not all modules are realizable: for example if $\G = C_{p^k}$
for $k \geq 3$, then the module $\bZ[C_p]\oplus \bZ[C_{p^i}]$,
$2\leq i < k$,
is not realizable by any smooth $\G$--action on $\nCP$.
However, we have the following ``stable" realization result.
We say that two subgroups $\pi_1$ and $\pi_2$ of $\pi$ are
\emph{disjoint} if $\pi_1\cap\pi_2 =\{1\}$.
\begin{thmc}
 Let $\cS$ denote a set of subgroups
of $\G=C_m$ containing at most two maximal elements (under inclusion
of subgroups). If $\cS$ has two maximal elements, suppose
that they are disjoint. Then
there exists an integer $N = N(\cS)\leq m$ such that
any permutation module 
$$ \bigoplus\{ \bZ[\G/\H_\alpha]^{k_\alpha}: \H_\alpha \in \cS\}$$
is realizable by a smooth $\G$--action on some $\nCP$, provided
that the multiplicity $k_\alpha > N$ for each of the maximal stabilizer subgroups. 
\end{thmc}

\begin{ack}  Our interest in permutation actions on homology
 was stimulated by a talk of A Edmonds \cite{Ed1} at the CRM, Montr\'eal (August  2000) on the
 orbit types of locally linear  topological actions on positive definite
 $4$--manifolds. This
research was partially supported by NSERC Discovery Grant A4000. The first author also wishes to thank the SFB 478, Universit\"at M\"unster, for hospitality and support.
\end{ack}

\section{The linear models}
 We investigate the orbit structure and possible permutation actions
 that can
arise on an equivariant connected sum of linear $\pi$--actions on
$\cp$. Let $(X, \pi)$ denote such a connected sum (always assuming that $X$ is simply-connected), and fix a generator $t$ for
the cyclic group $\pi =C_m$ of odd order $m$. 

\subsection{Linear actions on $\cp$}
We have already
pointed out that a linear $\pi$--action on $\cp$ (in standard form) has three fixed points $p_1 = [1, 0, 0]$, $p_2 = [0, 1, 0]$, and
$p_3 =[0, 0, 1]$. At these points, the tangential isotropy 
representations are described by the rotation numbers
$(a,b)$, $(-a, b-a)$ and $(-b, a-b)$, considered as 
pairs of integers mod $m$. For example, $(a,b)$ stands for
the linear $\pi$ representation $\bC(\zeta^a) \oplus \bC(\zeta^b)$.
The condition for an effective action  of $\pi$ on $\cp$ is
just that the greatest common divisor $(a,b,m) =1$. 
Up to oriented
equivalence, we can't distinguish the representations 
$(a,b)$ from $(b,a)$ or $(-a,-b)$ so these pairs of rotation numbers
are all identified. It will be convenient to use the notation
$\cp(a,b;m)$ for this linear action, or just $\cp(a,b)$ if the
cyclic group is understood.

Next we observe that the standard linear actions on $\cp$ always have
at least three invariant complex lines (topologically $2$--spheres), namely
the spans $S_1 = [p_2,p_3]$, $S_2 = [p_1,p_3]$ and 
$S_3 = [p_1,p_2]$. These invariant $2$--spheres are contained in 
the fixed sets of the
subgroups of order $m_1 = (a-b,m)$, $m_2 = (b,m)$ and 
$m_3= (a,m)$ respectively. Notice that the orders $m_1$, $m_2$
and $m_3$ are all co-prime, so that distinct $2$--spheres in the singular
set  of $(\cp, \pi)$ are fixed by distinct subgroups of $\pi$,
which intersect only in the identity element.
Any three divisors of $m$ can be obtained this way 
by appropriate choice of the rotation numbers, so we can obtain
actions with up to five distinct orbit types. Apart from the  fixed points and free orbits, the other three isotropy subgroups may be chosen 
arbitrarily, subject only to the condition that any two intersect in the identity. In the special case when $m = p^k$ for some prime $p$,
the divisibility condition implies that there are at most three distinct
orbit types. For any $\cp(a,b)$ the fixed set of each subgroup
of $\pi$ consists either of (i) three isolated points ($p_1$, $p_2$, and $p_3$),
or (ii) one isolated point and a disjoint $2$--sphere (one of the pairs
$p_i$, $S_i$).

\subsection{Linear actions on $S^4$} We take $S^4$ as
the unit sphere in $\bR^5 = \bR \oplus \bC \oplus \bC$ and write
the action of our generator in standard diagonal form:
$$ t\ \mapsto\ \left ( \vcenter{\xymatrix@R-28pt@C-24pt{
1&&\cr&\zeta^a&\cr&&\zeta^b}}\right )$$
as a matrix in $SO(5)$.
The action is again effective if and only if $(a,b,m) =1$.
We will use the notation $S^4(a,b; m)$ for this action.

For these actions on $S^4$, 
we always have two fixed points $p_1 = (1,0,0)$
and $p_2(-1,0,0)$ with the rotation numbers $(a,b)$ and $(a,-b)$.
In addition, the unit $3$--sphere $S^3 \subset 0 \oplus \bC \oplus \bC$
 always has two $\pi$--invariant circles from the two complex
factors, so we get two invariant $2$--spheres in the action.
These are $S_1 = \{(x, z_1,0)\vv x^2 + |z_1| = 1\}$
and $S_2 =  \{(x, 0, z_2)\vv x^2 + |z_2| = 1\}$, which are just the fixed
sets of the subgroups of order $m_1 = (a,m)$ and $m_2 = (b,m)$.
Again these two isotropy subgroups intersect only in the identity,
and we can obtain up to four orbit types in general.
If $m = p^k$ for some prime $p$, then three orbit types is the 
maximum possible. 
For any $S^4(a,b)$ the fixed set of each subgroup
of $\pi$ consists either of (i) two isolated points ($p_1$ and $p_2$),
or (ii)  a  $2$--sphere (either $S_1$ or $S_2$).

One common feature of these two linear models can be seen in any
smooth action on a connected sum.
\begin{theorem}\label{thm: fixedsetstruct}
If a cyclic group $\pi$ of odd order acts smoothly on a closed, smooth
simply-connected $4$--manifold $X$ with positive definite intersection
form, then the fixed point set $\Fix(X, \pi')$ for any non-trivial
subgroup $\pi'\subset \pi$ is a union of
isolated points and $2$--spheres.
\end{theorem}
\begin{proof} Since $X \simeq \nCP$ for some $n \geq 0$, it follows
that there is a standard basis $\{e_1, \dots, e_n\}$ for $H_2(X;\bZ)$
on which $\pi$ acts by permutations. Let $\pi'\neq 1$ be a non-trivial
subgroup of $\pi$, and let $C_p \subseteq \pi'$ be a subgroup of odd
prime order $p$. Since the action is orientation-preserving, $\Fix(X,
C_p)$ is a disjoint union of isolated points and oriented
surfaces. Since $\pi$ and therefore the subgroup $C_p$ acts by
permutations, the decomposition of $H_2(X;\bZ)$ as a $\bZ[C_p]$--module
has no summands of cyclotomic type (see \cite[Proposition 1.1]{Ed} for this
terminology). Therefore, by \cite[Proposition 2.4]{Ed} the fixed set has
zero first homology, and the result follows.
\end{proof}

\subsection{Equivariant connected sums}
To define the \emph{equivariant connected sum} of two smooth 
$4$--dimensional 
$\pi$--manif\-olds $(X, \pi)$ and $(Y, \pi)$, 
we select fixed points, $x \in X$ and $y \in Y$ with
rotation numbers of the form $(a, b)$ on $T_x X$ and
 $(a, -b)$ on  $T_y Y$. This just means that the tangent
representations are equivalent by an
\emph{orientation-reversing} isomorphism. We now
construct as usual the connected sum by removing small
$\pi$--invariant disks centered at the fixed points, and
equivariantly identifying concentric
annuli around $x$ and $y$ via the exponential map.
 Then $X \shp Y$ becomes a
smooth $\pi$--manifold in the standard way.

More generally, if $(Y,\pi')$ is a smooth $\pi'$--action,
for some subgroup $\pi'\subseteq \pi$, we will
describe the connected sum of $(X, \pi)$ with $(Z,\pi)=(\pi \times_{\pi'}Y)$ along an orbit of $\pi'$--fixed points. 
Select a point $x\in
X_{(\pi')} := \{x\in X \vv \pi_x = \pi'\}$,  and a point
$y \in \Fix(Y, \pi')$ whose rotation numbers (mod $|\pi'|$)
agree with those of $x$ up to orientation-reversal as
above. Then we perform the connected sum operation
$\pi/\pi'$--equivariantly, matching the points of the
orbit of $x$ in $X$ with the orbit $\pi/\pi'\times \{y\}$.
The resulting smooth $\pi$--manifold is denoted
$X \shp (\pi\times_{\pi'} Y)$.

We will refer to either of these operations as the equivariant
connected sum. The first is just a special case of the second
when $\pi' = \pi$. Notice that when $x \in X$ lies on a
$2$--sphere fixed by $\pi'$, then the rotation numbers at
$x$ have the form $(a,0)$ mod $|\pi'|$. The orientation-reversing
matching condition implies that the rotation numbers at
$y \in Y$ equal $(-a,0) \equiv (a,0)$ mod $|\pi'|$, so the point
$y$ lies on a fixed $2$--sphere in $Y$.

Since the $\pi'$--rotation numbers are constant along
connected components of $\pi'$--fixed sets,  we have some
choice in selecting the points at which to perform 
the connected sum.
\begin{lemma} The $\pi$--diffeomorphism type of
the equivariant connected sum
$X \shp (\pi\times_{\pi'} Y)$ is independent of the 
choices of $x$ and $y$ within connected components
of $X_{(\pi')}$ and $\Fix(Y, \pi')$.
\end{lemma}
\begin{proof}
This follows from an equivariant version of the isotopy extension
theorem (see \cite[VI.3]{Br}). 
\end{proof}

\begin{remark}
In an equivariant connected sum of linear actions on $\cp$
(or $S^4$),  the points of
$\pi/\pi'\times \{x\}$ and $\pi/\pi'\times \{y\}$ in the
construction described above will be called \emph{connecting
points}. If one connects up copies of $(\cp, \pi)$ using only
$\pi$--fixed connecting points, the induced action on the homology
of the connected sum obtained is trivial. Otherwise there exists
at least one non-trivial permutation module $\bZ[\pi/\pi']$ in
the homology of the connected sum.
\end{remark}

\begin{definition}
Let $\pi'\subseteq \pi$ be a subgroup of order $m'\mid m$.
We say that the linear
actions $X =\cp(a,b;m)$ and $Y=\cp(a',b';m')$ have a \emph{pair
of matching $\pi'$--fixed components} if there exists components
$C_1 \subset X_{(\pi')}$ and $C_2 \subset \Fix(Y,\pi')$ whose
rotation numbers on $T_x X$ at $x \in C_1$ and on
$T_y Y$ at $y\in C_2$ agree (mod $|m'|$) after an
orientation-reversing equivalence.
\end{definition}

\subsection{Trees, roots and branches}
We are going to associate a \emph{weighted tree}
to each equivariant connected sum of linear actions on $\cp$.
Given a positive odd integer $m=|\pi |$ and a positive integer $n$, we will consider $\pi$--equivariant trees $\bT$ with $n$ (type I), or $n+1$ vertices (type II),
having the following properties (see \cite[\S 2]{Ser1}):
\renewcommand{\labelenumi}{{\rm(\roman{enumi})}}
\begin{enumerate}
\item The vertices and edges of $\bT$ are permuted by the $\pi$--action,
preserving the incidence relation.
\item There is a $\pi$--fixed vertex called the \emph{root} vertex
for  the tree. We number the vertices $V = \{v_1, \dots, v_n\}$ (type I)
or $V = \{v_0, v_1, \dots, v_n\}$ (type II) so that the vertex $v_1$ (respectively
$v_0$ in type II) is the root vertex.
\item For a type II tree, the root vertex $v_0$ is the \emph{unique}
$\pi$--fixed vertex.
\item There is a $\pi$--invariant partial ordering on the vertex set $V$ 
($v_i < v_j $ implies $g\cdot v_i < g\cdot v_j$ for all $g \in \pi$), and the
root is the unique minimal element in this partial ordering.
\item The edge set $E$ is $\pi$--invariantly \emph{directed} so that 
for each $e\in E$, the initial and terminal vertices $v_i = \bd_0 e$ and $v_j =\bd_1 e$ satisfy $v_i < v_j$.
\end{enumerate}

We are going to think of the vertices $\{v_1, \dots, v_n\}$ of $\bT$  as representing
standard linear actions $\cp(a_i,b_i;m_i)$, for $1\leq i \leq n$. 
The root vertex $v_1$ of
a type I tree will have $m_1 =m$. For the type II trees,
the root vertex $v_0$ will be of the form $S^4(a_0,b_0;m_0)$,
with $m_0 =m$. We will translate from divisors of $m$  to 
subgroups by letting $\pi_i$ denote the subgroup of order $m_i$.
This discussion should motivate the following definition.

\begin{definition} An \emph{admissible,  weighted}
tree $(\bT, \pi)$ is a $\pi$--equivariant, part\-ially-ordered, directed tree with the properties listed above, and in addition:
\begin{enumerate}
\item Each vertex $v_i$ has \emph{weights} $(a_i,b_i;m_i)$,
where $m_i$ is a divisor of $m$, and $a_i$, $b_i$ are a
 pair of integers whose common divisor $(a_i, b_i, m_i)=1$.
\item If $v_i < v_j$, then $m_j$ divides $m_i$.
\item The weights of $v_i$ are the same as those of $g\cdot v_i$,
for all $g\in \pi$.
\item If $v_i = \bd_0 e$ and $v_j =\bd_1 e$ for some edge $e \in E$,
then the linear actions $\cp(a_i,b_i;m_i)$ and $\cp(a_j,b_j;m_j)$ have a \emph{pair of matching $\pi_j$--fixed components} 
\end{enumerate}
\end{definition}

We use the same equivalence relation on the weights as for
rotation numbers. In other words, in the triple $(a_i, b_i;m_i)$
the pair $(a_i,b_i)$ is well-defined only mod $m_i$, and
$(a_i,b_i) = (b_i,a_i) = (-a_i,-b_i)$. Let $\cR$
denote the collection of weights for the vertices of $(\bT, \pi)$,
and let $\bbT = (\bT, \pi, \cR)$ denote an admissible weighted
tree as defined above. There is an obvious notion of
\emph{equivalence} between weighted trees $\bbT$ and $\bbT'$,
involving a  bijection of vertex sets and edge sets which is compatible with the $\pi$--action, incidence relations,  and the weights. The ordering structure
is useful internally for describing the construction of
equivariant connected sums from this data, but one could attach
the various $\cp$'s is many different orders.
\begin{lemma} Given an admissible, weighted tree $\bbT$,
the set $\bT_0$ of $\pi$--fixed vertices and edges is
an  admissible, weighted subtree $\bbT_0\subseteq \bbT$.
\end{lemma}
\begin{proof}
In particular, we are asserting that $\bT_0$ is \emph{connected},
and of course it contains the root vertex by definition.
The proof is immediate from the partial ordering property: any
 $\pi$--fixed vertex $v_i$ is connected to the
root vertex by a path of edges. It follows that $m_j =m$ for
any intermediate vertex on the path, so the path lies
in $\bT_0$.
\end{proof}
We will call this $\pi$--invariant subtree the \emph{homologically trivial}
subtree of $\bbT$, since it gives rise to a homologically
trivial action. Notice that the complement $\bbT \setminus \bbT_0$
is a disjoint union of $\pi$--orbits of
admissible, weighted trees for various
subgroups  $\pi'\subset \pi$. We call these the \emph{branches}
in $\bbT_0$. For each branch, we can repeat the process to
find the homologically trivial $\pi'$--subtree of the branch. This gives
a canonical way to decompose the tree into simpler pieces.

\begin{theorem} Given an admissible, weighted tree $\bbT$,
there is an equivariant connected sum of linear actions
$(X(\bbT), \pi)$ such that  $X(\bbT) = \nCP$, up to
diffeomorphism. Conversely,
given an equivariant connected sum $(X, \pi)$ of linear actions
with $X = \nCP$, there
is an admissible, weighted tree $\bbT$ such that $(X(\bbT), \pi)
= (X,\pi)$. Two trees $\bbT$ and $\bbT'$ give $\pi$--equivariantly
diffeomorphic actions if and only if they are equivalent.
\end{theorem}
\begin{proof} Given an admissible, weighted tree $\bbT$
we can construct an action on $\nCP$, using the data given by the weights. We call  these  $(X(\bbT), \pi)$ \emph{tree manifolds}.
The last condition
on edges allows us to connect the vertex $v_i = \cp(a_i,b_i;m_i)$
to $v_j = \cp(a_j,b_j;m_j)$ by  equivariant connected sum.  To see
that equivalent trees give rise to equivariantly
diffeomorphic actions, we start with the homologically
trivial subtree $\bbT_0$. By the ``matching" condition on
edges, we see that the weights at every vertex in $\bbT_0$
are determined by the weights at any single vertex. This is because
we are taking the connected sum at $\pi$--fixed points, whose
rotation numbers within each $\cp(a,b;m)$ are in the set
$(a,b)$, $(-a,b-a)$, and $(-b, a-b)$. In forming the equivariant
connected sum we pick one of these, say $(a,b)$ and attach
$\cp(a,-b;m)$. Since $\bbT_0$ is connected, we determine
all the rotation numbers in $\bbT_0$ by this process. 

Conversely, suppose that $(X, \pi)$ is an equivariant connected
sum of actions on $\nCP$. We will argue by induction on $n$, 
where the case $n=1$ is clear. We could also encounter the
case $n=0$, which is just $S^4(a,b;m)$. By construction, there exists
a $\pi$--orbit $\{g\cdot e\vv g\in \pi\}$ in $H_2(X;\bZ)$, where
$e\in H_2(X;\bZ)$ is represented by some $\bC P^1\subset \cp$
used in the equivariant connected sum operation. This factor
is $\pi'$--invariant, where $\pi'\subset\pi$
is the stabilizer of $e$. We may choose $\pi'$ a minimal element
among the set of isotropy subgroups. Let $(X_0, \pi)$
be the result of removing a $\pi$--orbit of a $\pi'$--equivariant
tubular neighbourhood $\nu(\bC P^1\subset \cp)$ from $X$,
and then gluing in a $\pi$--orbit of $\pi'$--invariant $4$--disks. We
obtain $(X_0, \pi)$, which by induction is a disjoint union of tree-manifolds and $(Y,\pi) =\pi\times_{\pi'} \cp(a',b';m')$.
Since $(X,\pi) = (X_0,\pi)\shp (Y,\pi)$, it remains to see that the
connecting points $\pi/\pi'\times \{ x\}$ in $X_0$ may be
$\pi$--equivariantly isotoped into one of the tree factors. But
this is clear, since $\Fix(X_0,\pi')$ is a disjoint union of
isolated point and $2$--spheres. If $\pi'=1$ there is nothing
to prove since we can move the orbit of $x$ around in the free
part of the action. If $\pi' \neq 1$, then the orbit of $x$ either
lies in disjoint branches of the tree, or may be isotoped into
a $\pi'$--fixed $2$--sphere in one of the vertices. 
This shows that $(X,\pi)$ is a tree manifold.
\end{proof}

\subsection{Isotropy and fixed $2$--spheres}
We will now discuss the isotropy subgroups and the existence
of invariant $2$--spheres for the tree manifolds. 
A \emph{$\pi'$--isotropy $2$--sphere} is a two-dimensional
component of $X_{(\pi')}$.
Suppose first
that $\bbT_0$ is a homologically trivial tree with root
$\cp(a,b)$, and let 
$X= X(\bbT_0)$. We may describe the rotation numbers at 
the other vertices algebraically by considering the three possible
choices $p_1$, $p_2$, $p_3$ for attaching edges. With the
conventions above, the possible new weights are $w_1=(a,-b)$,
$w_2=(-a,a-b)$ and $w_3=(-b, b-a)$. We can define the matrices
$$K = \mmatrix{1}{\hp 0}{0}{-1}, 
\qquad L = \mmatrix{-1}{\hp 0}{\hp 1}{-1}, \qquad 
R=\mmatrix{\hp 0}{-1}{-1}{\hp 1}$$
with the property that $w_1 = Kw$, $w_2 = Lw$ and $w_3=Rw$,
where $w = (a,b)$ and we perform the matrix multiplication
with the weights as column vectors. The following result is clear.

\begin{lemma} Let $\Gamma  \subset GL_2(\bZ)$
be the subgroup generated by the matrices $K$, $L$ and $R$. If $\bbT_0$
is a homologically trivial tree with root weights $w=(a,b)$,
then the weights at the other vertices are of the form
$\gamma\cdot w$ mod $m$ for some $\gamma \in \Gamma$.
\end{lemma}
In principle, we could now determine all the possible isotropy
subgroups in $X(\bbT_0)$ from this recipe.
Notice the relations
$$L^k = (-1)^k \mmatrix{\hp 1}{0}{-k}{1}, \qquad 
R^k = \mmatrix{\hp f_{k-1}}{-f_k}{-f_k}{\hp f_{k+1}}$$
where $f_k$ denotes the $k^{th}$ Fibonacci number starting with
$f_0 =0$ and $f_1 =1$. Here is one simple observation.
\begin{lemma}
There is a homologically trivial tree $\bbT_0$ such that
$X(\bbT_0)$ contains a $\pi'$--isotropy $2$--sphere for every subgroup
$\pi'$ of $\pi$.
\end{lemma}
\begin{proof} Starting with the weights $(1,0)$ and applying
$L^k$ for $1\leq k \leq m$ we will obtain the weights 
 $w_i = (1,-m_i)$ for each divisor $m_i$ of $m$. Since the
linear action $\cp(1,-m_i;m)$ has a $\pi_i$--isotropy $2$--sphere
(where $|\pi_i| = m_i$) we are done.
\end{proof}
\begin{remark}
This result is of course very inefficient. We would like to prescribe
a collection of subgroups $\{\pi_1,\dots,\pi_r\}$ and find
the \emph{minimal} integer $n\geq 0$ such that there exists
a homologically trivial action on $\nCP$ with a $\pi'$--isotropy
$2$--sphere exactly when $\pi' = \pi_i$ for some $i$, $1\leq i\leq r$.
\end{remark}
For the general case, we can analyse the isotropy groups
and invariant $2$--spheres by considering the complement
$\bbT \setminus \bbT_0$ as a disjoint union of weighted trees.
The same process applies to each branch with respect to the
stabilizer subgroup $\pi'$ of the branch. To the operation
of the matrices $K$, $L$, and $R$ above we add the operation
of reduction mod $m'$ as we enter the branch. 
\begin{lemma} Let $X(\bbT)$ be a tree manifold. Then the weights
at the root vertex determine the weights at all the vertices.
\end{lemma}

\subsection{Permutation modules}
Let $X =X(\bbT)$ be a tree manifold and consider the
action of $\pi$ by permutations on the homology group
$$ H_2(X;\bZ) = \bigoplus\{ \bZ[\G/\H_\alpha]^{k_\alpha}: \H_\alpha \in \cS\}$$
We will assume that $X$ is an action on $\nCP$ for some $n>0$
(thus eliminating a type II tree with only the root vertex).
The set $\cS$ is the set of stabilizer subgroups for elements
in the standard basis $\{e_1,\dots,e_n\}$ of $H_2(X;\bZ)$.
We observe that if $\{1\} \in \cS$, then $\bZ[\pi]$ is obtained by
equivariant connected sum along a free orbit of points.
At the other extreme,  $\pi \in \cS$ if and only if the homologically
trivial part $\bbT_0\subset \bbT$ is non-empty (i.e. $\bT$ has
type I). The subgroup $\pi_\alpha =\pi$ contributes a summand
$\bZ^k$ to $H_2(X;\bZ)$, where the Euler characteristic
$\chi(\Fix(X,\pi)) = k+2$. This means that
in the equivariant connected sum $X(\bbT_0)$
there are exactly $k$ vertices.
\begin{lemma}
The proper, non-trivial, subgroups in the set $\cS$ are exactly the set of subgroups $\{\pi'\subset \pi\}$
for which (i) $\Fix(X,\pi')$ contains a $2$--sphere at a vertex $v_j$,
and (ii) there is an edge $e\in \bT$ with $\bd_0 e = v_i$
and $\bd_1 e = v_j$ such that $m_i >  m_j=|\pi'| $. 
\end{lemma}
\begin{proof}  Suppose
that $1\neq \pi' \in \cS$. Then there exists a vertex $v_j$ with
stabilizer $\pi' =\pi_j$ of order $m_j $. The conditions follow 
easily.
\end{proof}

We will be interested in the \emph{maximal} elements of
the set $\cS$ under the partial ordering by inclusion on subgroups.
Recall that disjoint subgroups are those which intersect only
at the identity element.
\begin{lemma} Let $X = X(\bbT)$ be a tree manifold diffeomorphic
to $\nCP$ for some $n>0$.
Then the set $\cS$ of stabilizer subgroups for the permutation
modules in $H_2(X;\bZ)$ has the unique maximal element
$\{\pi\}$ if $\bT$ has type I, and otherwise $\cS$ has at most
two maximal elements, and these are disjoint.
\end{lemma}
\begin{proof}
Of course if $\bT$ has type I, then the root vertex is $\pi$--fixed
and the homologically trivial part $\bbT_0$ is non-empty, so we
get copies of $\bZ$ in the permutation module.
If $\bT$ has type II, then there are no $\pi$--fixed vertices except
the root vertex, but this vertex is of the form $S^4(a,b;m)$
so its invariant $2$--spheres are homologically zero. However,
the equivariant connected sum must be performed at a $\pi$--orbit of connecting points, either in the free part of the action or in one
of the (at most) two singular $2$--spheres. The isotropy groups
of these singular $2$--spheres in $S^4$ are the maximal elements
in $\cS$, which are necessarily disjoint subgroups of $\pi$.
\end{proof}
\begin{remark}
Since the linear actions $\cp(a,b)$ have at most three non-trivial
proper isotropy groups, the realizable sets $\cS$ of
stabilizer subgroups have an additional property for type I 
tree manifolds with the root vertex as the \emph{unique} 
$\pi$--fixed vertex: $\cS$ contains at most three maximal \emph{proper}
subgroups, $\pi_1$, $\pi_2$ and $\pi_3$, and $\pi_i \cap \pi_j =
\{1\}$ for $i\neq j$. 
\end{remark}

We now have two stability results for realizing permutation modules
by tree manifolds.
\begin{theorem}\label{thm: permone}
 Let $\cS$ denote a set of subgroups
with $\pi \in \cS$.
There exists an integer $N = N(\cS)\leq m$, such that
any permutation module 
$$ \bigoplus\{ \bZ[\G/\H_\alpha]^{k_\alpha}: \H_\alpha \in \cS\}$$
with multiplicity $k_0 > N$ for the trivial representation 
$\bZ = \bZ[\pi/\pi]$ is realizable by a  tree manifold $X(\bbT)$.
\end{theorem}
\begin{proof}
We start with any homologically trivial action that contains
an isotropy $2$--sphere for each $\pi_\alpha \in \cS$. This can be
achieved by taking enough vertices in $\bbT_0$, as
shown above. We then use these $\pi_\alpha$--isotropy
two spheres as the attaching spheres for $k_\alpha$
branches of the form $\pi\times_{\pi_\alpha}\cp$.
\end{proof}
For type II trees we  assume stability at both maximal
elements of the set of stabilizer subgroups.
\begin{theorem}\label{thm: permtwo}
 Let $\cS$ denote a set of subgroups
with at most two maximal elements $\{\pi_1,\pi_2\} \in \cS$.
If $\cS$ has two maximal elements,
suppose that they are disjoint. Then 
 there exists an integer $N = N(\cS)\leq m$ such that
any permutation module 
$$ \bigoplus\{ \bZ[\G/\H_\alpha]^{k_\alpha}: \H_\alpha \in \cS\}$$
with multiplicities $k_i>N$, for $i=1,2$,
is realizable by a type II tree manifold $X(\bbT)$.
\end{theorem}
\begin{proof}
In the type II case we must attach a $\pi_i$--homologically
trivial branch to each of the $\pi_i$--isotropy $2$--spheres
in the root $S^4$. If we take sufficiently many vertices
in these branches we can obtain isotropy $2$--spheres for all the remaining subgroups of $\cS$. This is because each $\pi_\alpha \in \cS$
is contained in one of the maximal elements $\pi_1$ or $\pi_2$.
Then we attach orbits of the form $\pi_i\times_{\pi_\alpha} \cp$
in the $\pi_i$--branch. This realizes the module $\bZ[\pi_i/\pi_\alpha]$
in the homology of the branch, but the $\pi$--equivariance
gives us the module
$$\ind _{\pi_i}^{\pi}(\bZ[\pi_i/\pi_\alpha]) = 
\ind _{\pi_i}^{\pi}( \ind_{\pi_\alpha}^{\pi_i}(\bZ)) = 
\ind_{\pi_\alpha}^{\pi}(\bZ) = \bZ[\pi/\pi_\alpha]$$
in the homology of $X$.
\end{proof}
\begin{remark}
The method of proof shows that any $N \geq m$ works. More
precisely, we could take $N = \max\{ |\pi_\alpha| : \pi_\alpha \in \cS\}$.
These estimates are probably far from ``best possible".
\end{remark}

\begin{example}
The permutation module  $\bZ[C_{15}] \oplus \bZ[C_{21}]\oplus
\bZ[C_{35}]$, with $\pi=C_{105}$,
is \emph{not} realizable as $H_2(X;\bZ)$ for any equivariant
connected sum of linear actions. In this case
$\cS$ contains the three maximal elements $\{C_3,C_5,C_7\}$.
However, the module $\bZ \oplus \bZ[C_{15}] \oplus \bZ[C_{21}]\oplus
\bZ[C_{35}]$ is realizable, starting from the root $\cp(10,3;105)$.
\end{example}

\section{Equivariant moduli spaces}

In 
\cite{HL2} the theory of Yang--Mills moduli
spaces (\cite{Doa}, \cite{DK}) 
was adapted to the equivariant setting. 
We give an informal sketch of the main features and refer
to these sources for details.
Let $P\to X$ be a principal $SU(2)$--bundle over a
smooth, closed, oriented, simply-connected $4$--manifold $X$. 
Let $\cA(P)$ denote the space of $SU(2)$ connections
on $P$, and $\cB(P) = \cA(P)/\cG(P)$ the quotient space of
connections by the action of the gauge group $\cG(P)$.
If we fix a Riemannian metric on $X$, we can decompose
the space of Lie algebra valued forms
$$\Omega^2(X;\ad P) = \Omega_{+}(X;\ad P)
\oplus \Omega_{-}(X;\ad P)$$
into eigenspaces of the $\ast$--operator. The curvature
operator 
$$F\colon \cA(P) \to \Omega^2(X;\ad P)$$
is gauge invariant and  decomposes as $F = F^+\oplus F^-$.
We can define the \emph{self-dual} Yang--Mills moduli space
$$\cM(P) = \{[A]\in \cB(P) \vv F^-_A = 0\}$$
Now suppose that we have a smooth, orientation-preserving
action $(X, \pi)$.
We fix a real analytic structure on $X$ compatible with
its smooth structure and the given $\pi$--action, and a
 real analytic $\pi$--invariant metric on $X$ (see \cite{IL}).
Since action of $\pi$ on $X$ preserves the orientation,
and the bundle $P$ is classified by 
$c_2(P) \in H^4(X;\bZ)$, for
each element $g\in \pi$, we can find a (generalized)
bundle map $\varphi_g\colon P \to P$ such that
the diagram
$$\xymatrix{P \ar[r]^{\varphi_g}\ar[d]& P\ar[d]\cr
X \ar[r]^g& X}$$
commutes. Let $\cG(\pi) :=\{\varphi_g \vv  g\in \pi\}$
denote the group of generalized gauge transformations.
Since the indeterminacy in the choice of $\varphi_g$ is
an element of $\cG$, we have a short exact sequence
$$1 \to \cG \to \cG(\pi) \to \pi \to 1$$
of groups, and a well-defined $\pi$--action on $\cB$.
The gauge-invariance of the curvature operator implies
that the map $F$ is also $\pi$--equivariant, and hence
we obtain a natural $\pi$--action $(\cM, \pi)$ on 
the Yang--Mills moduli space.

Without further modifications, the moduli space could
be a very singular object. The index $\delta_A$ of its fundamental
elliptic complex
\addtocounter{dummy}{1}
\begin{equation}\label{eqn: ellcx}
\xymatrix{\Omega^0(X;\ad P) \ar[r]^{d_A}&
\Omega^1(X;\ad P) \ar[r]^{d^-_A} &\Omega_{-}(X;\ad P)}
\end{equation}
is called the \emph{formal dimension} of the moduli space
at $[A] \in \cM$, but it need not be its geometric dimension.
The operator
$d_A$ is the covariant derivative associated to the 
connection $A$, and the operator $d^-_A$ is the linearization
of the curvature.

The homology groups $H^i_A$, $0\leq i\leq 2$, of the fundamental
elliptic complex are finite dimensional real vector spaces,
whose dimensions $h^i(A) = \dim H^i_A$ appear in the
formula 
$$\delta_A = h^1(A) - h^0(A) - h^2(A)\ .$$
By the Atiyah--Singer Index Theorem,
the formal dimension has a purely topological
expression
$\delta_A = -8c_2(P) -3(1 - b_2^-)$, where $b_2^-
= \dim H^2_{-}(X;\bbR)$. In particular, it is independent of the choice of
base point $[A]$.  A necessary condition for $P$ to admit
any self-dual connections is $c_2(P)\leq 0$, and we will
only use the case $c_2 = -1$ and $b_2^- =0$ in this paper.

In the equivariant setting, the elliptic complex inherits
an action of the stabilizer $\cG_A(\pi)\leq \cG$ of the connection
under the generalized gauge group action. 
The stabilizer
sits in an exact sequence
$$1\to \Gamma_A \to \cG_A(\pi) \to \pi_A\to 1$$
where $\pi_A$, by definition, is the image of $\cG_A(\pi)$
in $\pi$. Here $\Gamma_A$ is the stabilizer of $A$
in $\cG$, and there are just two possibilities:
$\Gamma_A = \{\pm 1\}$ if $A$ is \emph{irreducible}, or
$\Gamma_A = S^1$ if  $A$ is \emph{reducible}. The
latter holds when the structural group $SU(2)$ of 
$P$ reduces to $S^1$, or equivalently, when the associated
complex vector bundle $E \to X$ splits as $E = L \oplus L^{-1}$
for some complex line bundle over $X$.

The Kuranishi method gives finite dimensional ``local charts"
$$\phi_A\colon H^1_A \to H^2_A$$
at each connection $[A] \in \cM$, where $\phi_A$ is a 
$\cG_A(\pi)$--equivariant smooth map with $(d\phi_A)_0 =0$.
Then in a neighbourhood of $[A]$, the moduli space $\cM$
is locally isomorphic to $\phi_A^{-1}(0)/\Gamma_A$.

\subsection{Equivariant general position}
\label{thm: modulifeature1}
In \cite{HL2},  an equivariant perturbation of the 
Yang--Mills equations was constructed (based on  the method of Bierstone 
\cite{Bie}) to obtain an equivariant
 \emph{general position} moduli space of 
self-dual connections $(\cM,\pi))$. After perturbation,
the equivariant moduli space is locally described by
finite dimensional charts $\phi_A\colon H^1_A \to H^2_A$ as above, which
are in Bierstone general position with respect to the $\pi_A$ action
(if $A$ is irreducible) or the  $S^1 \times \pi_A$ action (if $A$
is reducible). If $H^2_A =0$ after the perturbation, then the moduli
space is  \emph{equivariantly transverse} at $[A]$. We obtain the following
good properties:

\begin{list}{(\roman{enumi})}{\usecounter{enumi}
\setlength{\leftmargin}{.8truecm}}

\item $(\cM,\pi)$ is an equivariantly Whitney stratified space
\cite{M}  with an effective $\pi$--action and open smooth 
manifold strata
$$\cM^*_{(\pi')} = \{ x\in\cM^*\vv  \pi_x = \pi' \subseteq 
\pi\},$$
where $\cM^* \subseteq \cM$ is the subset of irreducible 
connections (up to gauge equivalence).

\item For  $\pi' \subseteq \pi$ each component of the 
fixed point set $\Fix(\cM^*,{\pi'})$ is the moduli space 
of $\pi'$--invariant connections on $P$ with respect to a 
$\pi'$--$SU(2)$ bundle structure on $P$ 
(compare \cite{FS}, \cite{FL}).

\item The strata have topologically locally trivial 
equivariant cone bundle neighbourhoods in $(\cM,\pi)$.

\item  $(\cM,\pi)$ has an equivariant compactification 
$(\bbar\cM,\pi)$. When $c_2(P) = -1$, and $X$ has positive
definite intersection form, 
$$(\bbar\cM,\pi)= (\cM\cup \tau(X\times [0,\lambda_0)),\pi)$$ where the ideal boundary $\bd\bbar\cM =\tau(X\times 0)$ is  the Taubes embedding of $X$ as the set of ideal 
``highly-concentrated" connections. Under this embedding,  the boundary of $\bbar\cM$ has  a smooth equivariant collar 
neighbourhood, which is $\pi$--equivariantly diffeomorphic to $X \times [0,\lambda_0)$ with the product action ($0<\lambda_0
\ll 1$).

\item The dimensions of the strata  $\cM_{(\pi')}$ can be 
computed as the index of the $\pi'$--fixed set of the fundamental elliptic complex.  A stratum $\cM_{(\pi')}$ is non-empty
whenver its  formal dimension is positive. In particular,
when $c_2(P) = -1$, the \emph{free} stratum $\cM^*_{(e)}$
is a smooth, non-compact $5$--dimensional manifold.
We may assume that the free stratum is connected.
\end{list}
\begin{remark}
For the case $c_2=-1$ and 
$X$ has positive definite intersection form, the moduli
space is equivariantly transverse on the subspace
$\cM_{\lambda_0}:=\cM^*\,\cap\,\tau(X\times [0,\lambda_0))$.
of highly concentrated connections \cite[\S 9]{FU}. 
\end{remark}

From now on we assume that $X$ is
a closed, smooth, simply-connected $4$--manifold with
positive definite intersection form. In other words, we assume that
$X\simeq \nCP$ for some $n>0$. We also restrict attention to
the equivariant moduli space $(\cM,\pi)$ for $c_2 = -1$,
 and let $\pi=C_m$ for $m$ odd. Under these assumptions, the
moduli space has a number of additional properties. Most of
these properties were established in \cite{HL1} for the special
case when the induced $\pi$--action on homology is trivial. We
will adapt  the statements to remove this assumption, and indicate
where the proofs need to be generalized.

\subsection{Reducible connections}
For $X \simeq \nCP$ we fix a standard basis $\{e_1,\dots, e_n\}$
for $H_2(X;\bZ)$ so that $Q_X(e_i,e_j) = \delta_{ij}$.
\begin{lemma}{\rm\cite[2.2.6]{DK}}\qua
Let $P \to X$ be a principal $SU(2)$--bundle with $c_2(P) = -1$,
where $X \simeq \nCP$.
There is a $\pi$--equivariant bijection $e_i \leftrightarrow [D_i]$
between the set $\{e_1, e_2, \dots, e_n\}$  of  standard
basis  elements for $H_2(X)$ and the set
$\cR :=\{[D_1], [D_2], \dots, [D_n] \}$ of gauge equivalence classes of
reducible $SU(2)$--connections on $P$. 
\end{lemma}
\begin{proof} The reducible connections correspond
to splittings $E = L \oplus L^{-1}$, for some complex line
bundle $L  \to X$. Such splitting  exist if and only if $c_2(E) = 
-x^2$ for some $x \in H^2(X;\bZ)$, since any such class  
determines a line bundle $L$ with $c_1(L) =x$. Since the only
solutions for the equation $x^2 =1$ (up to $\pm 1$) come from the Poincar\'e
duals $\hat e_i$ of the basis elements $e_i$, we have a bijection
between the two sets. However, we have
a $\pi$--action on the set of reducibles and on the homology
classes, and we must check that the actions correspond.
Let $\pi_i$ denote the stabilizer of $e_i$ under the $\pi$--action
on homology. Then, by the same proof as \cite[Corollary 1]{HL1},
there exists a $\pi_i$--equivariant $S^1$ bundle $L_i \to X$
such that $c_1(L_i) = \hat e_i$. Since an $S^1$--bundle has
a unique self-dual connection $[D_i]$, up to gauge
equivalence, we see that
$[D_i]$ is $\pi_i$--invariant. This shows that $[g^*D_i] = [D_j]$ if
and only if $c_1(g^*L_i) = L_j$, which holds
if and only if $g^*(e_i) = e_j$. Therefore the actions agree.
\end{proof}
The most striking feature of this $c_2=-1$ equivariant moduli space
is the existence of equivariantly  transversal local charts at each
of the reducible connections. Here we generalize \cite[Theorem 15]{HL1}.
\begin{theorem}\label{thm: eqcones}
Let $P \to X$ be a principal $SU(2)$--bundle with $c_2(P) = -1$,
where $X \simeq \nCP$. Then for each reducible connection $[D]$,
there is a $\pi_D$--invariant neighbourhood $\cN_D$ of $[D]$ in $(\cM,\pi)$ on which the moduli space is equivariantly transverse.
The stratified space $(\cN_D, [D])$ is  $\pi_D$--equivariantly homeomorphic to the cone over
some linear action of $\pi_D$ on $\cp$. Away from the cone point,
these two stratified spaces are equivariantly diffeomorphic.  
\end{theorem}
\begin{proof} The proof of this result in the homologically trivial
case in \cite[\S 4]{HL1} is fairly complicated, but the changes needed
to allow a non-trivial action on homology are minimal. The point is that
the set of reducibles is finite, so that a $\pi$--invariant neighbourhood
of a reducible $[D]$ has the form $\cU = \pi\times_{\pi_D} \cN$
for some $\pi_D$--invariant neighbourhood $\cN$ of $[D]$.
We may assume that $\cN\cap g(\cN)=\emptyset$  if $g\notin \pi_D$,
so the quotient group $\pi/\pi_D$  permutes
a set of disjoint copies of $\cN$ freely and transitively. The perturbations needed to
achieve $\pi_D$--equivariant transversality
 in a neighbourhood of $[D]$
(exactly as carried out in \cite[\S 4]{HL1}) can be extended
$\pi$--equivariantly to achieve transversality at each reducible
in the $\pi$--orbit of $[D]$. Since the set of reducibles is a disjoint
union of $\pi$--orbits, and the perturbations are local, we can
achieve transversality at each orbit separately.
\end{proof}
\begin{remark} We have shown that $H^2_D =0$ for each
reducible connection $[D]$. The local cone $\cN_D$
is just the quotient of $H^1_D = \bC(\chi_1)\oplus
\bC (\chi_2) \oplus \bC (\chi_3)$ by the stabilizer
$\Gamma_D = S^1$ acting diagonally as complex multiplication. 
The weights $\chi_i$, $1\leq i\leq 3$ are linear characters of $\pi_D$,
and the subspaces $\bC(\chi_i)$ are just the fixed
sets $\Fix(H^1_D, \pi_D(\chi_i))$, where 
$\pi_D(\chi_i):=\{(g, \chi_i(g^{-1})\vv g\in \pi_D\} \subset \pi_D\times S^1$. We call the boundary $\lk(D) :=\bd\cN_D$ the \emph{link}
of the reducible connection.
\end{remark}

\subsection{Orientation} The Yang--Mills moduli space
inherits an orientation from that of $X$ (see  \cite{Dob}, \cite[\S 5.4]{DK}).
Donaldson shows that the real determinant line bundle 
$\Lambda(P)$ associated to the elliptic complex
$$\OmegaP\colon \quad 0 \to \Omega^0(\ad P) \to \Omega^1(\ad P) \to 
\Omega_{-}(\ad P)\to 0$$
has a canonical trivialization over $\cB$ which induces the 
given orientation on $X$ times the inward pointing normal,
where $X$ is embedded as the 
Taubes boundary in $\cM$. We call this the canonical
orientation of $\cM$.
\begin{lemma}[see {\cite[Lemma 8]{HL1}}]\label{lemma: orientability}
Let  $\pi=C_m$ for $m$ odd, and for any subgroup
$\pi'$ of $\pi$
let $\cC \subset \cM^*_{(\pi')}$ be a connected component.
Then the  canonical orientation on $\cM$
induces a preferred orientation on the smooth manifold $\cC$.
\end{lemma}
\begin{proof}
If $[A]\in \cC$, we can split the elliptic complex
$$\OmegaP = 
\OmegaP^{\pi'} \oplus \bigl [\OmegaP^{\pi'}\bigr ]^\perp$$
into a fixed subcomplex and a perpendicular complex.
It follows that the line bundle $\Lambda(P) = 
\Lambda_{t(\cC)}\otimes \Lambda_{n(\cC)}$, where $\Lambda_{t(\cC)}$ 
is the determinant line bundle of $\OmegaP^{\pi'}$ and 
$\Lambda_{n(\cC)}$ is for the complementary part of $\OmegaP$.
Since $\pi$ is odd order cyclic, the action of $\pi'$ on
$ \bigl [\OmegaP^{\pi'}\bigr ]^\perp$ induces a complex structure, 
and hence a preferred orientation on $\Lambda_{n(\cC)}$.
Then the canonical orientation on $\cM$ induces an orientation on $\Lambda_{t(\cC)}$ as well. However, 
 the moduli space is locally modelled
on a  the zero set of a smooth map $f\colon H_A^1 \to H_A^2$
 in Bierstone  general position.
Furthermore, the fixed set
$(\Coker df_0)^{\pi'}=0$ and $\Ker df_0$ (which is fixed 
under $\pi'$) is the tangent space to the manifold stratum $\cC$ 
at $[A]$. Therefore $C$ has a preferred orientation.
\end{proof}
We may apply this to the free stratum, and consider the
copies of $\cp$ which bound $\pi_D$--invariant neighbourhoods
$\cN_D$ of each of the reducible connections $[D]\in \cM$.
Note that we have an $S^1$ determinant line bundle
$\mathcal L_{\pi'}$ over $\Fix(\cM,\pi')$ equal to
 the top exterior power with respect to the complex structure
on 
$ \bigl [\OmegaP^{\pi'}\bigr ]^\perp$.
\begin{lemma}{\rm\cite[Example 4.3]{Dob}}\qua
The induced orientation on each copy of $\cp =\bd\cN_D$ linking
a reducible connection $[D]$ is the standard complex orientation.
\end{lemma}
\begin{proof}
 This amounts
to the statement that the preferred complex line bundle 
$\mathcal L_{\pi_D}$, restricted to the link of each 
reducible fixed by $\pi_D$, is just the Hopf bundle over $\cp$.
\end{proof}
\begin{corollary} \label{cor: loop}
Let $\pi_D$ be the isotropy group of a reducible
connection $[D]$, and let $\pi'\subseteq \pi_D$.
\begin{enumerate}
\item
Then there is no continuous 
path $\gamma\colon [0,1]\to
 \wbar{\cM_{(\pi')}}$ such that $\gamma(0) = \gamma(1) = [D]$,
with $\gamma(0,1) \,\cap\, \cN_D$  a disjoint union of
one-dimensional strata.
\item Suppose that $[D'] = g^*[D]\neq [D]$ for some $g\in \pi$. Then
there is no  continuous 
path $\gamma\colon [0,1]\to
\wbar{\cM_{(\pi')}}$ such that $\gamma(0) = [D]$, $\gamma(1) =[D']$,
with $\gamma(0,1) \,\cap\, \cN_D$  and $\gamma(0,1) \,\cap\, \cN_{D'}$  a disjoint union of one-dimensional strata.
\end{enumerate}
\end{corollary}
\begin{proof} 
The argument is given in the proof of \cite[Theorem C, page 729]{HL1}, 
based on
the orientability of the $\pi_D$--fixed strata in the moduli space (see Lemma \ref{lemma:
orientability}). For part (i), suppose that there exists a such a closed
path $\gamma$ in $\Fix(\cM, \pi_D)$, intersecting the cone
$\cN_D$ in two  distinct one-dimensional strata with
isotropy groups $\pi_D(\chi_i)$, $\pi_D(\chi_j)$ respectively, in the 
$\pi_D\times S^1$ action of the local model.
But the $S^1$ determinant line bundle $\mathcal L_{\pi_D}$
restricted to these two strata is  positively oriented  on each stratum
by the complex structure. However the $S^1$ determinant
 line bundle  extends over
$\Fix(\cM, \pi_D)$, and so over $\gamma$ giving a
contradiction.
This  shows that
such a loop $\gamma$ can't exist.

For part (ii) we use the orientation-preserving 
action of $g\in \pi$ to identify the 
$S^1$ determinant line bundles $\mathcal L_{\pi_D}$ and
$\mathcal L_{\pi_{D'}}$ at the distinct reducibles $[D]$ and
$[D']$. Now the assumption that the $S^1$ determinant
 line bundle  extends over the path $\gamma$ again gives a
contradiction.
\end{proof} 

The existence of topologically locally trivial cone bundle
neighbourhoods for the strata in $\cM$ allows the
possibility of transporting rotation number information
from $X$ to the links of the reducible connections.
If $x\in X$ lies in $\Fix(X,\pi')$ then $(T_x X,\pi')$ splits into
a trivial $\pi'$--representation tangent to the fixed set,
and a complementary \emph{normal isotropy representation} $V_x$. If
$x\in \Fix(X,\pi')$ is an isolated fixed point, then $V_x$
is just the  isotropy representation $(T_x X,\pi')$ whose
eigenvalues are the rotation numbers (mod $|\pi'|$) at $x$.
Let $\cM_{\lambda_0} = \cM^*\cap \tau(X \times [0,\lambda_0))$
denote the intersection of the Taubes collar with
the subspace of irreducible connections. We will study
the  strata containing a singular set of the form
$F\times (0, \lambda_0)\subset\cM_{\lambda_0}$,
where $F \subset \Sing(X,\pi)$.

\begin{lemma}{\rm\cite[Corollary 3.2]{HL2}}\qua
\label{lemma: equivTaubes} Let  $F\subset X_{(\pi')}$ be
a connected component with normal isotropy
 representation $V_x \subset T_x X$  at $x \in F$.
 Then there exists a smooth manifold stratum  $\cC\subset
\cM^*_{(\pi')}$, such that
$\cC\, \cap\,\cM_{\lambda_0} = F \times (0, \lambda_0)$.
In addition, the
cone bundle neighbourhood for $\cC$ is a $\pi'$--equivariant disk
bundle with normal isotropy representation $V_x$. 
\end{lemma}
It follows that such a stratum carries a preferred 
orientation for its normal disk bundle in $\cM^*$. This can
sometimes be used to determine its closure. 
\begin{lemma}\label{lemma: arcorbit}
Suppose that $\pi'\neq \pi$ and that $x$ is an isolated fixed point in $\Fix(X,\pi')$. Then the closure  of the connected
$1$--dimensional stratum $\cC$ in $\cM^*_{(\pi')}$ 
containing $\{x\}\times (0,\lambda_0)$ either ends at
another isolated fixed point $y\neq x$ in $\Fix(X,\pi')$,
or ends at a reducible connection.
\end{lemma}
\begin{proof}
Under the given assumptions,  $\cC$ is a smooth non-compact 
$1$--dimen\-sional manifold with one  limiting endpoint at $x$. 
If its other limiting
endpoint $y\in X$ lies in the Taubes boundary, 
then $y$ must clearly be isolated and distinct from $x$. 
Suppose that $[A]\in \cM^*$
is the other limiting
endpoint in $\bbar\cC$. Then,  by general position, 
$[A]$ must have  a larger isotropy subgroup  $\pi_A\neq \pi'$.
 Therefore the full orbit of $\cC$
under the $\pi_A$ action must also have $[A]$ as a limiting  endpoint,
and for any $g \in \pi_A$, $g\notin \pi'$, the union
$\bbar\cC \,\cup\, g\cdot \bbar\cC$ gives a continuous path 
from $x$ to $gx \neq x$.
However, by transporting the preferred orientation
for $T_x X$ along this path, we  conclude that
the rotation numbers at $gx \in X$ are $(a,-b)$, opposite to
 those at $x$. This is a contradiction, since the $\pi$--action
on $X$ is orientation-preserving,   so we have eliminated
the possibility that $\cC$ has a limiting endpoint in $\cM^*$.
But $\cC$ is a smooth $1$--manifold, so  it has two ends,  
and $\bbar\cC$ must have another limit point in $\bbar\cM$. 
It follows that $\bbar\cC$ contains a reducible connection
\end{proof}

\begin{corollary}{\rm\cite[Corollary 3.3]{HL2}}\qua\label{cor: arcfromX}
Suppose that $x$ is an isolated fixed point in $\Fix(X,\pi')$ with
rotation numbers $(a,b)$. If the cancelling pair $(a,-b)$
of rotation
 numbers  does not occur in $(X,\pi)$,
then the closure  of the connected
$1$--dimensional stratum $\cC$ in $\cM^*_{(\pi')}$ 
containing $\{x\}\times (0,\lambda_0)$ ends at
 a reducible connection.
\end{corollary}
\begin{proof} If $\pi'=\pi$ the smooth $1$--manifold 
$\cC$ can't have a limiting endpoint in $\cM^*$, by general position,  so it must
give a smooth path to a reducible connection.
For the case $\pi'\neq\pi$, if the closure
of $\cC$ doesn't contain a reducible, then
 Lemma \ref{lemma: arcorbit}, asserts that $\bbar\cC$
gives a smooth path from $x\in X$ to another isolated
$\pi'$--fixed point $y\in X$. 
But then Lemma \ref{lemma: equivTaubes}
shows that the rotation numbers at $y$ are $(a,-b)$, 
which
is a contradiction.
\end{proof}
\begin{remark}\label{rmk: uniquearc} 
A smooth path $\gamma\colon [0, 1] \to
\bbar{\cM}$ emerging from an isolated
$\pi'$--fixed point on $X$, and ending at a reducible,
 is \emph{unique} (up to re-parametrization)
since for $0<t<1$, the path $\gamma(t)$
 is  contained in a $1$--dimensional smooth
manifold component of the stratum $\cM^*_{(\pi')}$.
We point out another situation where such a path must exist.
Recall that $\cR$ denotes the set of reducible connections.
\end{remark}
\begin{corollary}{\rm\cite[Lemma 17]{HL1}}\qua\label{cor: noreturn}
If $\Fix(\cR, {\pi'}) \neq \emptyset$, for some $\pi'\neq 1$ and $x\in X$
is an isolated $\pi'$--fixed point, then the closure of the
$1$--dimensional stratum in $\cM^*_{(\pi')}$ 
containing $\{x\}\times (0,\lambda_0)$
ends at a reducible connection.
\end{corollary}
\begin{proof} We remark that is possible for  all the reducible
connections to have have trivial stabilizer. For example, we
may start with $S^4(a,b)$ and form the equivariant connected
sum
with copies of $\cp$ along a free orbit. Then there exists a
 fixed arc in $\bbar\cM^*$ whose endpoints are two
isolated fixed points on $X$. 

Suppose, if possible, that there exists a $\pi'$--fixed smooth path
$\gamma$  joining two isolated $\pi'$--fixed points 
$x_1$, $x_2$ in $X$. Then the rotation numbers at $x_1$
and $x_2$ form a cancelling pair $(a,b)$ and $(a,-b)$. By 
Lemma \ref{lemma: equivTaubes}, there is a neighbourhood of
$\gamma$ in $\bbar\cM$ which is cut out  $\pi$--equivariantly
transversely. We may now pick a non-trivial subgroup
$C_p\subset \pi'$ and perturb the equations defining the moduli
space into $C_p$--general position. The perturbation can be
chosen to be the identity on this neighbourhood, so the arc
$\gamma$ does not change.
Since $\Fix(\cR, {\pi'}) \neq \emptyset$, 
we have a non-zero trace for the $C_p$--action
on $H_2(X;\bZ)$ and therefore $\chi(\Fix(X,C_p))\geq 3$.
It follows that there exists $x_3 \in \Fix(X,C_p)$ with
$x_3 \neq x_1, x_2$. Now the argument proceeds exactly
as in \cite[page 727]{HL1}.
\end{proof}

\subsection{Fixed sets in the closure of the free stratum}
For our applications to group actions on $\nCP$ we only
need to study the closure $\Me$ of the free
stratum of the moduli space. By general position
(see  \S \ref{thm: modulifeature1}), 
and the good cone structure at the reducibles (Theorem 
\ref{thm: eqcones}) the closures of any singular
strata of  dimension $\geq 5$ are disjoint from
the closure of the free stratum. This means that 
$\Me$ is the union of the free stratum
with (i) the reducibles and (ii) some singular strata of
dimension one or three. We are interested in
knowing more about the
closures of the singular strata inside $\Me$.

 \begin{lemma} {\rm\cite[Lemma 12]{HL1}}\qua
\label{lemma: conncomp}
Let $C\subseteq \cM^*_{(\pi')}$ be a non-compact,
connected component with $\pi'\neq 1$.
If $C$ is non-empty, then $\dim C < 5$. 
 If $\dim C = 3$, then the
closure $\bbar{C} \subset \bbar{\cM}$ must
intersect the Taubes boundary $\partial \bbar{\cM} = X
\times \{0\} \subset \bbar{\cM}$.
\end{lemma}
\begin{proof}
For  $\dim C \geq 5$, the proof given 
in \cite[Lemma 12]{HL1} applies without change: assuming
$C$ is non-empty leads to a contradiction.
If $\dim C =3$ and $\pi'\neq \pi$  we assume, if possible, that the closure $\wbar C$ is formed by adjoining some one-dimensional
strata in $\cM^*$ (with larger isotropy), together with some
$\pi'$--fixed reducible 
 connections. Now we  remove the intersection of
$\wbar C$ with the interiors
of the cones $\cN_D$ around each of these reducibles.
 We then obtain 
 a compact $\pi'$--fixed set
$(\cC, \bd \cC)$ in $\cM^*$, bounding a collection $F$
of $2$--dimensional $\pi'$--fixed
sets in the links. Suppose now that 
we  also have
\begin{enumerate}
\item another fixed set $(\cC', \bd\cC')$ in $\bbar\cM\setminus \cR$ bounding
another  collection $F'$ of $2$--dimensional $\pi'$--fixed
sets in the Taubes boundary $X$, and
\item the intersection $\cC\,\cap\, \cC'$ contains a singular stratum with larger isotropy group $\pi''\neq \pi'$.
\end{enumerate}
 In that case, following the outline of the argument given in 
\cite{HL1}, the next step is to pick a subgroup $C_p\subseteq \pi'$,
for some prime $p$,
and perturb into $C_p$--general position. 
The difficulty is that
the union $\cC\, \cup \,\cC'$
could now be perturbed into a smooth $C_p$--fixed
cobordism between the 
two collections of $2$--spheres, which is allowable.

However, if this ``bad" case occurs, the existence of such
a stratum $\cC'$ would imply that the sum of the two-dimensional
homology classes in $F'$ would map to zero under the
Taubes map $\tau_*\colon H_2(X;\bZ) \to H_2(\cM^*;\bZ)$.
Furthermore, $\mu([F']) = 0 \in H^2(\cB^*;\bZ)$ by 
\cite[5.1.2]{DK}. 
But the composite 
$$\xymatrix{H_2(X;\bZ) \ar[r]^{\mu}& H^2(\cB^*;\bZ)\ar[r]^{\tau^*}
&H^2(X;\bZ)}$$
is just the Poincar\'e duality map \cite[5.3.3]{DK}, so
$\tau^*(\mu([F']))=0$ implies that $[F']=0\in H_2(X;\bZ)$.
An easy modification of \cite[Corollary 4]{HL1} shows that,  whenever 
$\Fix(H_2(X;\bZ/p), C_p)\neq 0$, then 
 the  $2$--dimensional components 
of $\Fix(X,C_p)$ represent linearly independent elements of 
$H_2(X;\bZ/p)$. 
 It follows that $\tau_*([F'])\neq 0$
in $H_2(\cM^*;\bZ)$, and hence no such stratum $\cC'$
bounding $F'$ exists. 

Therefore, after the perturbation into 
$C_p$--general position, the perturbed\break bounding set $\cC$
will be a smooth $C_p$--fixed null-bordism for $F$ in $\cM^*$,
disjoint from the Taubes boundary.
The existence of such a null-bordism is now ruled out by
the $\mu$--map argument exactly as given  in \cite[Lemma 12]{HL1}.
\end{proof}

The following useful statements were proved in the argument just given.

\begin{corollary}{\rm\cite[Theorem 16]{HL1}}\qua
\label{prop: bounding}
\noindent
\begin{enumerate}
\item 
No non-empty collection of $2$--dimensional $\pi'$--fixed
sets in the links of reducible connections bounds a compact fixed
set in $\cM^*$.
\item If $\Fix(\cR, \pi') \neq \emptyset$, then no
 non-empty collection $F$ of $2$--dimensional $\pi'$--fixed
sets in $X$, with $[F]\neq 0$,  bounds a compact fixed set in $\cM^*$.
\item If there exists a $\pi'$--fixed $2$--sphere in $X$, representing
a non-zero homology class, then $\Fix(\cR,\pi')\neq \emptyset$.
\end{enumerate}
\end{corollary} 
\begin{remark}
The case $X=S^4(1,0)$ shows that the first assumption
 in part (ii) is
necessary. In general, if $\Fix(\cR,\pi') = \emptyset$ then
the trace of the action of  a generator of $\pi'$ on $H_2(X;\bZ)$
is zero, and so $\chi(\Fix(X, \pi') = 2$. It follows that
$\Fix(X, \pi')$ consists of two isolated points, or a single
null-homologous
$2$--sphere, which must bound a compact fixed set in $\cM^*$.
We will see in the next section that the assumption
$[F]\neq 0$ in part (ii) is not actually necessary.
\end{remark}

\begin{definition}
Let $[D]$ be a $\pi'$--fixed reducible connection, with
$\pi' \neq 1$. A 
$\pi'$--\emph{incident stratum} at $[D]$ is a connected
 component
$\cC$ of $\cM^*_{(\pi')}$ such that $[D] \in
\bbar{ \cC}$. The \emph{stabilizer} $\pi_{\cC}$
of $\cC$ is the subgroup of $\pi$ leaving the stratum
invariant.
\end{definition} 

A $\pi'$--incident stratum intersected with the cone
$\cN_D$ of the $\pi'$--linear action at  a reducible $[D]$
 is a connected component of
$(\cN_D\setminus [D])_{( \pi')}$. Since its closure has
 $[D]$ as a limit point,  $\pi'\subseteq\pi_D\subseteq \pi_{\cC}$.
\begin{remark}\label{rmk: Eulerlink}
There are at most three proper 
$\pi'$--incident strata at a reducible connection since $(\lk(D), \pi_D)$ is a linear action on
$\cp$. 
Their closures consist of  either (i) three $1$--dimensional $\pi'$--fixed strata intersecting  $\lk (D)$ in three isolated $\pi'$--fixed
points, or (ii)  a $1$--dimensional and a $3$--dimensional
$\pi'$--fixed stratum intersecting $\lk (D)$ in an isolated
 $\pi'$--fixed
 point and a $\pi'$--fixed $2$--sphere, respectively. 
It follows that the Euler characteristic 
$\chi(\Fix(\lk (D),\pi' ))=3$
for  each reducible $[D]\in \cR$, and each $\pi'\subseteq \pi_D$.
\end{remark}
\begin{proposition}\label{prop: stratatypes}
Let $\cC$ be a $\pi'$--incident stratum
with stabilizer $\pi_{\cC}$. Then
there exists a $\pi_{\cC}$--invariant neighbourhood $\nu(\bbar\cC)$
of $\bbar\cC\,\cap\,\cM^*$, such that
 $g(\nu(\bbar\cC))\,\cap\,\nu(\bbar\cC) =\emptyset$
for all $g\in \pi$, $g\notin \pi_{\cC}$. In addition, the disjoint union 
$\pi\times_{\pi_{\cC}} \nu(\bbar\cC)$ is a smooth $\pi$--invariant
submanifold of $\cM^*$ on which the moduli space
is equivariantly transverse.
The  closure $\bbar{\cC}$  is one of the following types:
\begin{enumerate}
\item $\bbar{\cC}\,\cap\, \bd\bbar\cM =\emptyset$, 
 $\dim \cC =1$, and $\cC$ has exactly two reducible limit points, 
\item $\bbar{\cC}\,\cap\, \bd\bbar\cM
\neq \emptyset$, $\dim \cC = 1$ or $3$,
and $\cC$ contains a unique  reducible limit point, or 
 \item $\bbar{\cC}\,\cap\,\bd\bbar\cM
\neq \emptyset$, $\dim \cC =3$ and
$\cC$ has more than one reducible limit point.
\end{enumerate}
\end{proposition} 
\begin{proof} We begin with the classification of incident strata into
types. 
First we consider  case (i) when
$\bbar{\cC}\,\cap\, \bd\bbar\cM =\emptyset$ and 
 $\dim \cC =1$. If $\pi'=\pi_D =\pi$, then such a $1$--dimensional
smooth stratum has no interior limiting endpoint in $\cM^*$ by
general position. In that case, its other limiting endpoint must
be another reducible $[D_j]$ (distinct from its initial limit
reducible $[D_i]$  by Corollary \ref{cor: loop}). 

If $\pi'=\pi_D\neq \pi$,
and if $\cC$ has an interior limiting endpoint $[A]\in \cM^*$,
then $\pi'\subsetneq \pi_A$ and
the orbit $\pi_A\times_{\pi_D}\cC$ is a $\pi_A$--invariant
collection of $1$--dimensional strata which are incident with the orbit $\pi\times_{\pi_D} [D]$
and have $[A]$ as a common limit point.
However, this possibility is eliminated by 
Corollary \ref{cor: loop}(ii).
It follows that our incident  $1$--dimensional stratum
$\cC$  has no interior limiting endpoint in $\cM^*$. Since we
also assumed that its closure didn't meet $\bd\bbar\cM$, the 
only remaining possibility is that its other limiting endpoint must
be a reducible $[D_j]\neq [D_i]$. That completes case (i).

Next we note that the possibility  $\bbar{\cC}\,\cap\, \bd\bbar\cM=\emptyset$ and $\dim \cC =3$ is ruled out by
  Lemma \ref{lemma: conncomp}.
In the remaining cases (ii) and (iii) we
suppose that $\bbar{\cC}\,\cap\, \bd\bbar\cM\neq \emptyset$.
When $\cC$ has a unique reducible limit point, both
$\dim \cC=1$ and $\dim \cC =3$ both occur in the linear models.
On the other hand, if $\cC$ has more than one reducible
limit point, and some other limit point in $X$, then it can't be $1$--dimensional, since a $1$--manifold
has at most two ends. It follows in case (iii) that $\dim\cC =3$.

Now that we have the classification of incident strata into types, 
we can observe that in each case $\cC$ is incident with
a part of the moduli space (either $\cN_D$ or $\cM_{\lambda_0}$)
where  the $\cM^*$ is cut out equivariantly transversely.
But then Lemma \ref{lemma: equivTaubes} implies the
existence of a  $\pi_{\cC}$--invariant 
neighbourhood $\nu(\bbar\cC)$ with the
required properties.
\end{proof} 
\begin{corollary}\label{cor: disjointstrata}
The closures of distinct $\pi'$--incident strata can intersect only
at reducibles.
\end{corollary}
\begin{proof}
For each limit reducible $[D]$, 
the group $\pi_D$ operates freely, away from the ``zero section"
$\bbar\cC$,
within each of the $\pi_D$--invariant neighbourhoods
$\nu(\bbar\cC) \subset \cM^*$.
 It follows that $\cC$ is the unique $\pi'$--fixed
stratum in $\nu(\bbar\cC)$.
\end{proof}


\section{The structure of the singular set}
The main result of this section is Corollary \ref{cor: connfixedset}, showing 
that there is a nice subset of  $\Fix(\bbar\cM, \pi')$ 
containing $\Fix(\cR,\pi')$ and $\Fix(X,\pi')$ which is
 \emph{path connected}
for each $\pi'\subseteq \pi$. From this we can deduce the
proof of Theorem B and prepare for the proof of Theorem
A in the next section.

\subsection{Connecting the fixed-point sets}
In the last section we studied the individual $\pi'$--incident strata at
the reducible connections. Let $\cC(\pi')$ denote the union
of the $\pi'$--incident strata at all $[D]\in \Fix(\cR,\pi')$, assuming
that $\Fix(\cR,\pi')\neq\emptyset$. Notice that the $\pi$--action
permuting the reducibles induces a $\pi$--action 
on $\cC(\pi')$ which permutes the type (i) strata. Each type (ii)
or (iii) stratum $\cC$ has a stabilizer subgroup $\pi_{\cC}$
containing $\pi'$.

For $\pi'\neq 1$, we now consider the following
union of $\pi'$--fixed sets:
$$\cF(\pi') :=\Fix(\cR,\pi') \,\cup\, \Fix(X,\pi') \,\cup\,
\bigcup\{\bbar\cC(\pi'')\vv \pi'\subseteq \pi''\}
$$
lying in the compactified moduli space $(\bbar\cM,\pi)$.
By Corollary \ref{cor: disjointstrata} this is a disjoint union
except for inclusion of strata, and possible common limit points among the $\pi'$--fixed
reducibles. The main result of this section is that $\cF(\pi')$ is connected.

\begin{remark}
If $\Fix(\cR,\pi')=\emptyset$, then we define $\cF(\pi')$
as the connected component of $\Fix(\bbar\cM^*,\pi')$
containing $\Fix(X,\pi')$. This exists since $\Fix(X,\pi')$
consists either of two isolated points or a null-homologous
$2$--sphere, and either case these bound in $\bbar\cM^*$.
\end{remark}

If we remove from $\cF(\pi')$ the open cones around
each of the reducibles in $\Fix(\cR,\pi')$, we obtain
 a disjoint union of smooth $\pi'$--fixed cobordisms between
the singular sets in the links $\{\lk(D)\vv [D] \in \Fix(\cR,\pi')\}$
and the singular set $\Fix(X,\pi')$ (see Proposition \ref{prop: stratatypes}). Note that all the components
of $\Fix(X,\pi')$ are connected to the reducibles by
$\pi'$--incident strata (Corollary \ref{cor: noreturn}
and Proposition \ref{prop: bounding}).
We wish to compare the total Euler characteristic of the stratum
cobordisms at the two ends.

Let $\chi(\Fix(X,\pi')) = r+2$, where $r$ is the number of
$\pi'$--fixed reducibles, and let $\ell$ denote the number of homologically
zero $2$--sphere components in $\Fix(X,\pi')$. 
We will  assign \emph{weight} $\chi(\cC) =2$ to each
type (i) stratum, and weight $\chi(\cC)=2(k-1)$ to
each type (iii) stratum with $k$ reducible limit points.
The type (ii) strata have weight zero.
The quantity
$$\chi:= \sum\{\chi(\cC)\vv \cC \in \cF(\pi')\}$$
is the total Euler characteristic used up in joining various
reducibles by closures of $\pi'$--incident strata.
Next we let $\tau(\cC) = 2(t-1) $, where $t$
denotes the number of $2$--sphere boundary components 
in $\Fix(X,\pi')$ in a  type (ii) or type (iii) stratum  $\cC$ (otherwise
let $\tau(\cC)=0$), and define the \emph{excess}
$$\tau:=  \sum\{\tau(\cC)\vv \cC \in \cF(\pi')\}$$
We have the basic relation
$$3r -\chi +\tau +2\ell = r+2$$
by comparing Euler characteristics, but there is actually
enough information to compute these quantities:
\begin{lemma}\label{lemma: excess} Let $1\neq \pi'\subseteq \pi$.
 If $\Fix(\cR,\pi') \neq \emptyset$, then the weight
$\chi = 2(r-1)$ and the excess $\tau = 0$.
In addition, there are no homologically trivial $2$--sphere
components in $\Fix(X, \pi')$.
\end{lemma}
\begin{proof}
Each linear action
on $\cp$ contributes 
Euler characteristic 3 to the $\pi'$--incident
strata, so we have an ``excess" of $3r-\chi$ which must
equal the contribution from $\Fix(X,\pi')$. 
We can divide the
$r$ reducibles in $\Fix(\cR,\pi')$ into $s$ disjoint subsets by
defining two reducibles to be \emph{equivalent} if they are both
limit points of the same  stratum $\cC$.
 Suppose there
are $k_i$ reducibles in the $i^{th}$ subset, for  $1\leq i \leq s$, and
$\sum k_i = r$. Each subset accounts for $2(k_i -1)$ units in $\chi$.
Then the equation
$$\chi = \sum 2(k_i -1) = 2r - 2s = 2(r-1) +\tau +2\ell$$
gives $2(1-s) = \tau +2\ell\geq 0$. Since $s\geq 1$ we are done.
\end{proof}
\begin{corollary}\label{cor: 2sphere}
  Suppose that $\Fix(\cR,\pi)=\emptyset$,
and $\pi'\neq 1$ is a  proper subgroup which is maximal
with respect to the property
$\Fix(\cR, \pi')\neq\emptyset$. 
Then there exists a unique $\pi'$--incident stratum of type (iii), and a unique 
$2$--sphere component in $\Fix(X,\pi')$.
\end{corollary}
\begin{proof} 
Suppose if possible that there are no $\pi'$--incident
strata of type (iii). Then the set $\cF(\pi')$ contains $r$ reducibles
connected in pairs by $(r-1)$ type (i) strata. But $\pi$ acts on 
$\cF(\pi')$ by permuting the reducibles, and since $\pi$ has
odd order, there would exist a $\pi$--fixed reducible, contrary to
our assumption. Now if $\cC$ is a $\pi'$--incident stratum of type (iii),
it must intersect $X$ in a $\pi'$--fixed $2$--sphere $F$. If $\pi_{\cC}
\neq \pi$, then the orbit $\pi\times_{\pi_{\cC}} \cC$ would contribute
$2|\pi/\pi_{\cC}|$ to the Euler characteristic of $\Fix(X, \pi_{\cC})$,
contrary to our calculation of $\chi$. Therefore $\pi_{\cC} =\pi$
and $F$ is $\pi$--invariant.

Since $\Fix(\cR,\pi) =\emptyset$, the fixed set
$\Fix(X,\pi)$ consists of two isolated points $x_0$, $x_1\in F$.
Now suppose that there is another $2$--sphere component $F'$
in $\Fix(X,\pi')$. By Proposition \ref{prop: bounding} it is contained
in a $\pi'$--incident type (ii) or  (iii) stratum $\cC'$. We conclude
by counting as above, that $\pi_{\cC'} = \pi$, and it follows that
$F\,\cap\,F' = \{x_0,x_1\}$. But this contradicts the effectiveness
of the $\pi$--action at $x_0$ (by consideration of the tangential
isotropy representation $T_{x_0}X$). 
\end{proof}
\begin{corollary} If $\Fix(\cR,\pi) =\emptyset$, there
 are at most two maximal, proper subgroups $\pi_1$, $\pi_2$,
such that $\Fix(\cR,\pi_i) \neq \emptyset$. If two such
subgroups exist, then $\pi_1 \cap\pi_2 =\{1\}$.
\end{corollary}
\begin{proof}
Suppose that $\pi_1$ and $\pi_2$ are maximal proper
subgroups with the given property. By the last result, they each
have a fixed $2$--sphere $F_1$, $F_2$ respectively,
and $F_1\cap F_2$ gives just the two isolated points 
$\{ x_0,x_1\}$
in $\Fix(X,\pi)$. In order to have an effective action at
$x_0$, we must have $\pi_1\cap \pi_2 =\{1\}$ and no
other such subgroup can exist in the linear model
for $(\bR^4, C_m)$.
\end{proof}

For a type (iii) stratum $\cC$, suppose that 
$\cR(\cC)=\{[D_{1}],\dots, [D_{k}]\}$ is the set of its 
 reducible limit points. The stabilizer group $\pi_{\cC}$ 
of $\cC$ acts by permuting the reducibles in $\cR(\cC)$. 
We  identify the following sub-types:

\begin{lemma}\label{lemma: type3}
Suppose that 
$1\neq \pi'\subset \pi$.
For each $\pi'$--incident type (iii)  stratum, the set
$\cR(\cC)$ of reducible limit points,
under
the permutation action of the stratum
 stabilizer group, is either
\begin{description}
\item[\rm(iii-a)]  fixed by $\pi_{\cC}$, so $\pi_D =\pi_{\cC}$
for all $[D]\in \cR(\cC)$, and $\pi'=\pi_{\cC}$, or 
\item[\rm(iii-b)] a disjoint union of free  $\pi_{\cC}/\pi'$--orbits, with
$\pi'=\pi_D$ for all $[D]\in \cR(\cC)$, and
$\pi'\subsetneq \pi_{\cC}$, or
\item[\rm(iii-c)] a disjoint union
of free $\pi_{\cC}/\pi'$--orbits, 
together with some
 $\pi_{\cC}$--fixed points $[D]$, and $\pi'\subsetneq\pi_D=\pi_{\cC}$.
\end{description}
\end{lemma}
\begin{proof}
Each $\pi'$--incident type (iii) stratum
 $\bbar\cC$ has a $2$--sphere boundary component
in  $\Fix(X,\pi')$ on which
the  stabilizer subgroup $\pi_{\cC}$ either
acts trivially,or semi-freely with two isolated $\pi_{\cC}$--fixed points. 
In the first case,  $\pi_{\cC}$ also acts trivially on the stratum
$\cC$ and fixes all the reducibles in $\cR(\cC)$.
In the second case, $\pi'\subsetneq \pi_{\cC}$, and the
two isolated $\pi_{\cC}$--fixed  points either lie in $\Fix(X, \pi)$ and bound
a $1$--dimensional $\pi$--fixed set in $\cM^*$, or lie in distinct
$1$--dimensional type (ii) strata (for $\pi_{\cC}$) which
both have $\pi_{\cC}$--fixed limit reducibles $[D]$. In that case,
these  limit reducibles are among the set of 
$\pi_{\cC}$--fixed point in $\cR(\cC)$. The remaining limit
reducibles are permuted freely by $\pi_{\cC}/\pi'$.
\end{proof}

\begin{corollary}  Suppose that $\Fix(\cR,\pi)=\emptyset$,
and $\Fix(\cR, \pi')\neq\emptyset$ for some $\pi'\neq 1$. 
Let $\pi''\supseteq \pi'$ be a maximal subgroup
such that $\Fix(\cR, \pi'')\neq\emptyset$.
If $\cC_0$ is the unique $\pi''$--incident stratum of type (iii),
then $\pi_{\cC_0} = \pi$, and the  set of  limit reducibles
$\cR(\cC_0)$ is a disjoint union of free $\pi/\pi''$--orbits.
\end{corollary}
This stratum $\cC_0$ will be called the \emph{maximal}
type (iii-b) stratum for $\pi'$. It is uniquely determined by any 
subgroup $\pi'\subseteq \pi''$ such that 
$\Fix(\cR, \pi')\neq\emptyset$.

\subsection{The singular set tree}
In order to show that the sets $\cF(\pi')$ are connected, we will
study the configuration of reducibles and incident strata more
abstractly. 
\begin{definition}\label{def: graph}
For each subgroup $1\neq \pi'\subseteq \pi$, 
we will associate a graph 
$$\Gamma(\pi'):=\Gamma(V, E)$$
 whose vertex set $V := \Fix(\cR,\pi')$, provided that 
$\Fix(\cR,\pi) \neq\emptyset$. 
If $\Fix(\cR ,\pi) =\emptyset$, we adjoin one more $\pi$--fixed vertex $v_0$, called the \emph{root} vertex, which is 
 common to all the graphs
$\Gamma(\pi')$, $\pi'\subseteq \pi$.

The edge set $E$ is determined by the classification of
the $\pi'$--incident strata  in $\cF(\pi')$
into types (see Proposition \ref{prop: stratatypes}). 
For a $\pi''$--incident stratum of type (iii-a) with
$\pi'\subseteq \pi''= \pi_{\cC}$, we pick
an ordering
$\{[D_1], \dots, [D_k]\}$ of the limit reducibles,
and extend the ordering $\pi$--equivariantly over the
orbit $\pi\times_{\pi_{\cC}}\cR(\cC)$.
For a  $\pi''$--incident stratum $\cC$ of type (iii-c)
we pick a $\pi_{\cC}$--fixed reducible $[D]$,
called the \emph{branch} vertex,
and extend the choice by $\pi$--equivariance over the orbit of $\cC$. 

\smallskip
\noindent
Two distinct vertices  $[D_i]$, $[D_j]$ are joined by an edge
if and only if:
\begin{enumerate}
\renewcommand{\labelenumi}{(\alph{enumi})}
\item[(e-1)] there is
 a $\pi''$--incident stratum of \emph{type (i)}, with
$\pi'\subseteq \pi''$, 
and limit reducibles $[D_i]$ and $[D_j]$, or
\item[(e-2)] there is
 a $\pi''$--incident stratum $\cC$ of \emph{type (iii-a)}, with
$\pi'\subseteq \pi''= \pi_{\cC}$, limit reducibles
 $[D_i]$, $[D_j]$, and $j = i+1$ in the chosen ordering of
$\cR(\cC)$, or
\item[(e-3)] there is
 a $\pi''$--incident stratum $\cC$ of \emph{type (iii-c)}, with
$\pi'\subseteq \pi''\subsetneq \pi_{\cC}$, 
   limit reducibles $[D_i]$, $[D_j]$, such that 
 $[D_i]$ is the branch vertex in $\cC$ and  $\pi_{D_j}=\pi''$.
\end{enumerate}
If $\Fix(\cR,\pi)=\emptyset$, the root vertex $v_0$ is joined
to a vertex $[D]$ if  and only if $[D]$ is a limit reducible
in the maximal type (iii-b) stratum $\cC_0$ for $\pi'$.

\end{definition}
\noindent
Recall that a connected graph with no circuits is called a
\emph{tree}.
\begin{corollary}\label{cor: fixedtree}
For each subgroup  $1\neq \pi'\subseteq\pi$, the graph $\Gamma(\pi')$
is a tree.
\end{corollary}
\begin{proof}
The graph  $\Gamma(\pi')$ contains no closed circuits, by
Corollary \ref{cor: loop}, so we must show that 
$\Gamma(\pi')$ is connected. Let's first consider
the case where $\Fix(\cR,\pi)$\break$\neq\emptyset$.
We have the basic relation
$\chi = 2(r-1)$ where $r$ is the number of $\pi'$--fixed
reducibles. In our  Euler characteristic count, each edge
$e\in E$ arising from a type (i)
stratum contributes $\chi(e)=2$,  and each type (iii) stratum
with $k$ reducible limit points contributes $\chi(e) = 2(k-1)$,
so we can interpret the quantity $\chi$ as twice the number of edges
$e\in E$. In other words, there   are
exactly $(r-1)$ edges in $\Gamma(\pi')$.  But a tree with
$k$ vertices has exactly $(k-1)$ edges, so if $\Gamma(\pi')$
were disconnected into $\ell$ trees of $k_i$ vertices, $1\leq k_i \leq \ell$,
we would get 
$$\chi = \sum 2(k_i-1) = 2r - 2\ell\geq 2(r-1)\ .$$
Therefore $\ell = 1$ and $\Gamma(\pi')$ is connected.

If $\Fix(\cR,\pi)=\emptyset$, but $\Fix(\cR,\pi')\neq\emptyset$,
let $\pi''$ be  a maximal subgroup with this
property. 
Then the $k$ limit reducibles of the maximal
$\pi'$--fixed stratum $\cC_0$ of type (iii-b) 
are permuted in a disjoint union of $\pi/\pi''$--orbits.
This stratum adds $k$ edges to the graph $\Gamma(\pi')$, 
instead of $(k-1)$, and the total number of edges
is given by $\frac{1}{2}(\chi + 2) = r$ since $\cC$ is incident
with a $\pi'$--fixed $2$--sphere in $X$.
It follows by counting as before that $\Gamma(\pi')$ is connected.
\end{proof}
\begin{corollary}\label{cor: connfixedset} For each $\pi'\subseteq\pi$,  the set
$\cF (\pi' )$ is path connected.
\end{corollary}
\begin{proof} We may assume that
$\Fix(\cR,\pi')\neq \emptyset$.
 By construction of the graph $\Gamma(\pi')$,
and Corollary \ref{cor: disjointstrata},
 it is clear that  
$\Gamma(V,E)$ is  connected if $\cF(\pi')$ is connected.
On the other hand,  all the components
of $\Fix(X,\pi')$ are connected to the reducibles by
$\pi'$--incident strata. Therefore
if $\Gamma(\pi')$ is connected then
$\cF(\pi')$ is connected.
\end{proof}

The internal structure of the tree $\Gamma(\pi')$ is clarified
by considering
 the \emph{stabilizer} $\pi_{\cF}$
of a subset $\cF(\pi')$.
We will use the well-known fact that
a finite group of odd order acting on a tree always fixes a vertex.
\begin{lemma}
For each set $\cF(\pi')$, $1\neq \pi'\subseteq\pi$, there exists $[D]\in\Fix(\cR,\pi')$
such that $\pi_D = \pi_{\cF}$, and $\pi_D$ is a maximal element
in the set of subgroups $\{\pi_{D'}\}$ for all
$[D'] \in \Fix(\cR,\pi')$. Moreover, 
$\cF(\pi_D) = \Fix(\cF(\pi'), \pi_D)$.
\end{lemma}
\begin{proof}
We note first that $\pi_{D'} \subseteq \pi_{\cF}$
for all $[D']\in \Fix(\cR,\pi')$ because $\cF(\pi')$ is connected.
On the other hand, the $\pi_{\cF}$--action on the tree $\Gamma(\pi')$
has a fixed vertex $[D] \in\Fix(\cR,\pi')$, so $\pi_{\cF}\subseteq
\pi'$. Therefore $\pi_{\cF} = \pi_D$ for this reducible.
Since $\pi'\subseteq \pi_{\cF}=\pi_D$, we get the obvious
inclusion $\cF(\pi_D) \subseteq \Fix(\cF(\pi'), \pi_D)$.
But if $\cC$ is a $\pi''$--incident stratum at a $\pi_D$--fixed
reducible, for $\pi'\subseteq \pi''$, then either it is $1$--dimensional and fixed by
$\pi_D$ as well, or it is $3$--dimensional and the $\pi_D$--fixed
set of its closure is in $\cF(\pi_D)$.
\end{proof}
\medskip

We now want to consider the full singular set of $(\cM,\pi)$, or at
least those components whose closures intersect either the
links of the reducibles or the Taubes collar. We define
$$\cF(X,\pi): = \bigcup\{\cF(\pi') \vv 1\neq \pi'\subseteq \pi\}$$
and define the associated graph
$$\Gamma(X,\pi) := \bigcup \{ \Gamma(\pi')\vv 1\neq \pi'\subseteq \pi\}$$
The main result is: 
\begin{theorem}
The graph $\Gamma(X,\pi)$ is a tree and the singular set
$\cF(X,\pi)$ is connected.
\end{theorem}
\begin{proof}
Suppose that $[D]$ and $[D']$ are reducible connections. If 
$\pi_{D} \subseteq \pi_{D'}$ then $[D']\in \cF(\pi_D)$ so we can
connect them by a path in $\cF(\pi_D)$, or equivalently by an edge path in $\Gamma(\pi_D)$. If $\Fix(\cR,\pi)\neq\emptyset$, this
shows that any connection $[D]$ can be connected
 to any $\pi$--fixed reducible $[D']$, and we are done.
If $\Fix(\cR,\pi)=\emptyset$ but $\pi_D\neq 1$ for some
reducible connection $[D]$, then we let $\cC_0$ denote the maximal
type (iii) stratum for $\pi_D$. By Lemma \ref{lemma: type3},
the stabilizer $\pi'=\pi_{D'}$ for any reducible $[D']\in \cR(\cC_0)$
is the maximal subgroup containing $\pi_D$ with respect to
the property that $\Fix(\cR,\pi')\neq\emptyset$. Then 
$[D']\in \cF(\pi_D)$ and $[D']$ is connected to the root vertex
in $\Gamma(\pi_D)$. Alternately, $[D']$ is connected to
$\Fix(X,\pi)= \{x_0,x_1\}$ inside $\bbar\cC_0$. Therefore
$\Gamma(X,\pi)$ and $\cF(X,\pi)$ are connected, and
$\Gamma(X,\pi)$ is a tree.
\end{proof}

\subsection{The proof of Theorem B}

After this preparation we can now show that a smooth
pseudo-free action of an odd order cyclic group on $X\simeq
\nCP$ must be semi-free. This answers a question of 
A Edmonds \cite{Ed1}, who pointed out that the hypothesis
of a \emph{smooth} action is necessary. Edmonds showed
\cite[Theorem 5.4]{Ed1}
that $C_{25}$ can act locally linearly and pseudo-freely
on $X = \nCP$, for $n=10$, inducing the representation
$\bZ[C_5]\oplus \bZ[C_5]$ on $H_2(X;\bZ)$. In particular,
by computing traces, we see that this action is not semi-free.

We begin with the following observation.
\begin{lemma} Let $(X,\pi)$ be a pseudo-free smooth action
of $\pi=C_m$, $m$ odd, on $X\simeq \nCP$.
Then  there are only type (i) or $1$--dimensional type (ii) 
strata in the moduli space
$(\cM^*,\pi)$.
\end{lemma}
\begin{proof}
We have already established in Proposition
\ref{prop: stratatypes} that the $\pi'$--incident strata
come in three types, and those of type (ii) or (iii) intersect
$\Fix(X,\pi')$ in isolated points or $2$--spheres. By Proposition
\ref{prop: bounding} and Corollary \ref{cor: 2sphere},
any $\pi'$--fixed $2$--sphere in $\cN_D$ for a reducible
$[D]$ would imply the existence of a singular $2$--sphere
in $(X,\pi)$. Since the given action is pseudo-free, this can't occur.
\end{proof}
The usefulness of this observation is shown by the following
two remarks.
\begin{lemma}
\label{lemma: stabilvertex}
Suppose that $\cC$ is a $\pi'$--incident stratum of type (i), 
with limit reducibles $[D]$ and $[D']$. Then
$\pi_D = \pi_{D'}$.
\end{lemma}
\begin{proof}
Since $\cC$ is $1$--dimensional, its intersection with
$\cN_D$ is fixed by $\pi_D$, hence its intersection with
$\cN_{D'}$ is also fixed by $\pi_D$. Therefore  $\pi_D\subseteq \pi_{D'}$.
Similarly, $\pi_{D'}\subseteq \pi_{D}$.
\end{proof}
\begin{corollary} 
$\pi_D = \pi_{D'}$ for all $[D], [D']\in \cR$ with non-trivial stabilizer.
\end{corollary}
We now finish the proof of Theorem B by using again the
well-known fact that
a finite group of odd order acting on a tree always fixes a vertex.
In our case, the group $\pi=C_m$ of odd order acts on the
tree $\Gamma(X,\pi)$, so there is a vertex $[D]\in\cR$ which is
fixed by $\pi$. But this means that $\pi\subseteq \pi_D$
so $\pi_D = \pi$ for all $[D] \in \cR$ with non-trivial stabilizer. 
In addition, we could have some free $\pi$--orbits of reducibles,
but these don't contribute to the singular set.

It follows that
every type (ii) stratum is also fixed by $\pi$, so the limit points of
these strata in $X$ consist entirely of isolated $\pi$--fixed points.
On the other hand, by Corollary \ref{cor: noreturn} any $\pi'$--singular
point in $X$ is a limit point of some type (ii) stratum. Therefore
the singular set of $(X,\pi)$ consists of isolated $\pi$--fixed points,
and the action is semi-free.

\section{Stratified cobordisms and the proof of Theorem A}
In this section we will show that the moduli space contains
an equivariant connected sum $(X(\bbT),\pi)$ of linear actions on $\cp$,
and provides an equivariant stratified cobordism between 
$(X(\bbT),\pi)$ and our given action $(X,\pi)$. By construction, the action
$(X(\bbT),\pi)$ will have the same permutation representation
on $H_2(X;\bZ)$ as the given action. 
The stratified
cobordism will have smooth strata, and equivariant vector
bundle neighbourhoods. This will allow us to compare the
isotropy groups and rotation numbers in $(X,\pi)$ with those
of $(X(\bbT),\pi)$.

We begin by describing how to realize edges in $\Gamma(X,\pi)$
by thickened paths in $\bbar\cM^*$.
Recall that the tree $\Gamma(X,\pi)$ may have some edges arising from type (iii) incident strata, listed as cases (e-2) and (e-3) in Definition
\ref{def: graph}. If $\cC$ is a $\pi'$--incident stratum of type (iii-c)
with $\pi'\subsetneq \pi_{\cC}$, then there exists a branch vertex
$[D]\in \cR(\cC)$ and each of the  limit reducibles 
appearing in free $\pi_{\cC}/\pi'$--orbits are joined
by an edge to the branch vertex. We realize this geometrically by
choosing a smooth embedded path $\gamma$ in $\bbar\cC$ 
from $[D]$ to some other reducible $[D']\in \cR(\cC)_{(\pi')}$. The interior
of this path will be chosen to lie in $\cM^*_{(\pi')}$, with the property
that $g\cdot \gamma$ is disjoint from $\gamma$ (except at $[D]$)
for all $1\neq g\in \pi_{\cC}$. Now let $\nu_1(\cC)$ denote the union
of small tubular neighbourhoods in $\cM^*$. around the paths $\{g\cdot\gamma \vv g\in \pi_{\cC}\}$. We extend by equivariance
to the orbit $\pi\times_{\pi_{\cC}} \nu_1(\cC)$.

For a type (iii-a) $\pi'$--incident stratum $\cC$ with $\pi'=\pi_{\cC}$ we have
chosen an ordering $$\{[D_1],\dots, [D_k]\}$$ of the limit reducibles.
We choose disjoint smooth embedded paths $\gamma_i$ in $\bbar\cC$,
from $[D_i]$ to $[D_{i+1}]$, for $1\leq i \leq k-1$, and thicken
as before to define $\nu_1(\cC)$.
We again  extend by equivariance
to the orbit $\pi\times_{\pi_{\cC}} \nu_1(\cC)$.

If $\Fix(\cR,\pi)=\emptyset$ we will also have at most two
maximal type (iii-b) strata, associated to subgroups $\pi_1$, $\pi_2$
with $\pi_1\cap \pi_2 =\{1\}$, which are maximal with respect to
the property $\Fix(\cC,\pi')\neq\emptyset$. In this case, we 
have $\Fix(X,\pi) = \{x_0,x_1\}$, so we can define $\nu(x_0$ to be
a small $\pi$--invariant $5$--ball in the  Taubes collar centered
at the fixed point $(x_0,\lambda_0/2)$. The boundary of this
$5$--ball is a linear action $S_0:=S^4(a,b)$ having invariant $2$--spheres
with isotropy $\pi_1$ and $\pi_2$.
Now if $\cC_0$ is a
maximal type (iii-b) stratum for $\pi_1$, its limit reducibles are
a disjoint union of free $\pi/\pi_1$--orbits. We choose disjoint smooth paths $\gamma_i$, $\pi$--equivariantly, 
 from a disjoint union of free $\pi/\pi_1$--orbits of points in $\Fix(S_0,\pi_1)$ to the limit reducibles in $\cC_0$. Then we
thicken these paths as above, and adjoin $\nu(x_0)$
  to define $\nu_1(\cC_0)$. 

Finally,  there may exist some free $\pi$--orbits of reducibles.
If $\pi\times [D]$ is such an orbit, let $\nu_1(D)$ denote
the union of a $\pi$--invariant collection of disjoint, thickened,
smooth paths from $\{g\cdot [D]: g\in \pi\}$ to a free $\pi$--orbit
of points in some $\bd\cN_{D'}$, if there exists a $\pi$--fixed
reducible $[D']$, or in $S_0$ if $\Fix(\cR,\pi)=\emptyset$.

\begin{definition} We define a subset of $\bbar\cM^*$ containing
all the reducible connections. 
Let $$\cN(X,\pi) :=\bigcup\{\cN_D: [D]\in \cR\}$$
and 
$$\nu(X,\pi) := \bigcup \{\nu(\cC): \cC \text{\ type (i) stratum}\}
\cup\{\nu_1(\cC): \cC \text{\ type (iii) stratum}\}, $$
and then define
$$\cD(X,\pi) := \cN(X,\pi) \,\cup\, \nu(X,\pi)\,\cup\, \bigcup\{\nu_1(D): \pi_D = \{1\}\}
$$
\end{definition}
\begin{theorem} The boundary $(Y,\pi):=\bd\cD(X,\pi)$ is 
an equivariant
connected sum of linear actions on $\cp$, and $H_2(Y;\bZ)
\cong H_2(X;\bZ)$ as permutation modules.
\end{theorem}
\begin{proof} The equivariant moduli space $(\cM,\pi)$
has given us an admissible,\break weight\-ed tree $\bbT$ based on $\Gamma(X, \pi)$. The weights are given by the linear actions
on $\cp$ in the links $\ell(D)$ of the reducible connections, 
and by construction the permutation action on $H_2(X;\bZ)$
is realized by the permutation action on the reducibles, which
are all contained in $\cD(X,\pi)$. The edges of $\bbT$
are given by the strata $\cC$ of type (i) together with the
paths $\gamma$ constructed in the definition of the subsets
$\nu_1(\cC)$ for the type (iii) strata.
\end{proof}
\begin{proof}[The proof of Theorem A]
Let $W:=W(X,\pi)$ denote the complement of the interior
of $\cD(X,\pi)$ in $(\bbar\cM^*, \pi)$. Then $W$ is a stratifed
$\pi$--equivariant cobordism between $(X,\pi)$ and the
equivariant connected sum $(Y,\pi)$. In addition, the tubular 
neigbourhoods of the singular
strata in $W$ are just the intersections $W\cap \nu(\cC)$
for all type (ii) or (iii) strata in $\cM^*$, together with
a thickened tube from $D(x_0)$ in $X$ to $S_0$
in the case $\Fix(\cR,\pi)=\emptyset$. It follows
that the isotropy structure of $(Y,\pi)$ is the same
as that of $(X,\pi)$. In addition, 
 these tubular neighbourhoods are the total spaces of equivariant
vector bundles over the strata in $W$, so the rotation numbers
at singular points in $(Y,\pi)$ match up with rotation numbers 
in $(X,\pi)$.
\end{proof}
\begin{proof}[The proof of Theorem C]
This follows from Theorem A and Theorem \ref{thm: permone}.
Combining Theorem A with Theorem \ref{thm: permtwo} gives
the stable realization theorem for permutation modules in the case
when $\Fix(\cR,\pi)=\emptyset$, or equivalently, in the case when
there are no trivial summands $\bZ$ in $H_2(X;\bZ)$.
\end{proof}
 

\end{document}